%% file: LodayStein1.tex
\documentclass[11pt, fleqn, leqno]{article}

\input LodayStein0.tex

\begin{document}
\numberwithin{equation}{section}
\numberwithin{enumi}{subsection}
%\swapnumbers
\setcounter{section}{-1}
\large \centerline {\Large Parametrized Braid Groups of Chevalley
Groups}
%\bigskip \bigskip \par \hfill{\tt [\today,
%in progress, }
%
%\hfill{\tt ls02dec16.tex not for distribution, needs lsaction0.tex]}
\bigskip \bigskip \par

\centerline {\bf Jean-Louis Loday and Michael R. Stein}
\bigskip \bigskip \par

\begin{abstract} We introduce the notion of a braid group
parametrized by a ring, which is defined by generators and relations
and based on the geometric idea of painted braids. We show that the
parametrized braid group is isomorphic to the semi-direct product of
the Steinberg group (of the ring) with the classical braid group. The
technical heart of the proof is the Pure Braid Lemma \ref{pblan},
which asserts that certain elements of the parametrized braid group
commute with the pure braid group.

More generally, we define, for any crystallographic root system, a
braid group and a parametrized braid group with parameters in a
commutative ring. The parametrized braid group is expected to be
isomorphic to the semi-direct product of the corresponding Steinberg
group with the braid group. The first part of the paper (described
above) treats the case of the root system $A_n$; in the second part,
we handle the root system {$D_n$}. Other cases will be treated in the
sequel \cite{LS2}.
\end{abstract}

\N
\section{Introduction}
Suppose that the strands of a braid are painted and that the paint
from a strand spills onto the strand beneath it, modifying the color
of the lower strand as in the picture below.

%figure 1
\setlength{\unitlength}{0.6cm}
\centerline{
\begin{picture}(6,6)
\put(0.9,5.2){$u$}
\put(2.9,0.5){$u$}
\put(0.3,0.5){$v+au$}
\put(2.9,5.2){$v$}
\put(3,2.8){$a=$ coefficient of spilling}
\put(1,5){\line(0,-1){1}}
\put(3,5){\line(0,-1){1}}
\put(1,2){\line(0,-1){1}}
\put(3,2){\line(0,-1){1}}
\put(1,4){\line(1,-1){2}}
\put(3,4){\line(-1,-1){0.8}}
\put(1,2){\line(1,1){0.8}}
\end{picture}
}
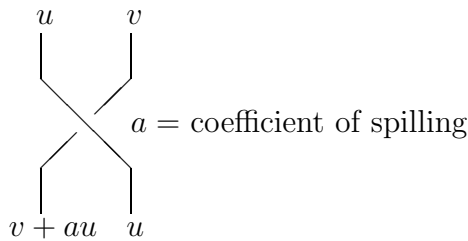
\centerline {figure 1}
\M

This gives rise, for any ring $A$, to the {\it parametrized braid
group} $Br_n(A)$, which is generated by elements $\y{i}{a}$, where $i$
is an integer, $1\leq i\leq n-1$, and $a$ is an element of $A$,
subject to the relations

\begin{align*}
(A1)\quad&& \y{i}{a} \y{i}{0} \y{i}{b}&= \y{i}{0} \y{i}{0} \y{i}{a+b}
& \\
(A1 \times A1)\quad&& \y{i}{a} \y{j}{b} &= \y{j}{b}
\y{i}{a}&\text{ if } |i-j|\geq 2,\\
(A2)\quad&&
\y{i}{a}\y{i+1}{b}\y{i}{c}&= \y{i+1}{c}\y{i}{b+ac}\, \y{i+1}{a} & \\
\text{ for any } a,b,c \in A \end{align*}

\N (A variation of this group first appeared in \cite{L}. The choice
of names for these relations will be explained in \S
\ref{explrelns} below. The derivation of these relations from the
painted braid model can be seen in Figures 2 and 3 of \S
\ref{explpix} below.)

Observe that when $A = \{0\}$ (the zero ring), one obtains the
classical Artin braid group $Br_n$, whose presentation is by
generators $\y{i}{}$, $1\leq i\leq n-1$, and relations

\begin{align*}
(A1 \times A1)\quad&& \y{i}{} \y{j}{}&= \y{j}{} \y{i}{}&\text{ if }
 |i-j|\geq 2,\\ (A2) \quad&& \y{i}{}\y{i+1}{}\y{i}{}&=
 \y{i+1}{}\y{i}{}\, \y{i+1}{} \\
\end{align*}

A question immediately comes to mind: does Figure 1 correctly reflect
the elements of the parametrized braid group? Up to equivalence, a
picture would be completely determined by a braid and a linear
transformation of the set of colors. This linear transformation lies
in the subgroup $E_n(A)$ of elementary matrices. Hence if the
parametrized braid group corresponds exactly to the pictures, it
should be the semi-direct product of $E_n(A)$ and $Br_n$. We will show
that this is almost the case: we need to replace $E_n(A)$ by the
Steinberg group $St_n(A) $ (\emph{cf.} \cite{St} \cite{Sb}).  \M

\N {\bf Theorem.} {\it For any ring $A$ there is an isomorphism
$$Br_n(A) \cong  St_n(A) \rtimes  Br_n\ ,$$
where the action of $Br_n$ is via the symmetric group ${\mathcal{S}_n}$.}
\M

The quotient of $Br_n(A)$ by the relation $\y{i}{0} \y{i}{0} =1$ is
the group studied by Kassel and Reutenauer \cite{KR} (in this quotient
group, our relation (A1) becomes $ \y{i}{a} \y{i}{0} \y{i}{b} =
\y{i}{a+b}$, which is exactly the relation used by \cite{KR} in place
of (A1)). They show that this quotient is naturally isomorphic to the
semi-direct product $St_n(A)\rtimes {\mathcal{S}_n}$ of the Steinberg
group with the symmetric group. So our theorem is a lifting of theirs.
\M

More generally, when $A$ is commutative, a Steinberg group $St(\Phi,
A)$ exists for any crystallographic root system $\Phi$. We construct a
{\it parametrized braid group}, $Br(\Phi, A)$, using a parametrized
version of the relations
$$ y_{\al} y_{\be} y_{\al}\cdots =  y_{\be} y_{\al} y_{\be}\cdots
$$
with $m(\al, \be)$ (see \S \ref{explm} for the definition)
generators on each side.  Our aim is to prove in all generality the
existence of an isomorphism
$$Br(\Phi,A) \cong St(\Phi,A) \rtimes Br(\Phi,0)\ .$$ Whenever such an
isomorphism exists, it implies an analogue of the theorem of Kassel
and Reutenauer for $\Phi$. 

\medskip
When $\Phi = A_n$, the only values of
$m(\al, \be)$ are $1,2$ and $3$, and the corresponding relations are
$(A1), (A1\times A1)$ and $(A2)$ respectively. The same is true for
the case $\Phi = D_n$ that we treat in \S \ref{Dn}.

More cases will be treated in the sequel \cite{LS2}.
\M

\N {\bf Content.} In \S \ref{explm} we outline the proof for a general
crystallographic root system $\Phi$. We define the parametrized braid
group $Br(\Phi, A)$ of type $\Phi$ with parameters in a commutative
ring $A$. For such data a Steinberg group $St(\Phi, A)$ has been
constructed in \cite{Sb} and \cite{St}. We explain how it is equipped
with an action of the Weyl group $W(\Phi)$, and (since $W(\Phi)$ is an
epimorphic image of $Br(\Phi)$), with an action of the braid group
$Br(\Phi) = Br(\Phi,0)$. We define a map $\phi : Br(\Phi, A)\to
St(\Phi, A)\rtimes Br(\Phi)$ which we conjecture is an isomorphism for
a general crystallographic root system $\Phi$, and we indicate a
procedure for proving this conjecture, which depends on a crucial
technical result called the Pure Braid Lemma (\ref{PBL}). Our proof
of the Pure Braid Lemma requires a manageable family of generators for
the pure braid group.  Such a family is easily obtained for the cases
$\Phi = A_{n-1}$ and $\Phi=D_n$ that we treat in this paper.

In \S \ref{explpix} we treat the case $\Phi= A_{n-1}$ (the
parametrized braid group presented in the introduction). Though this
case follows the general pattern outlined in the first section, \S
\ref{explpix} can be read independently. Moreover, in this case we
need not assume that $A$ is either commutative or unital.

In \S \ref{Dn} we treat the case $\Phi= D_n$, which is technically
more difficult due to the complexity of the generators of the pure
braid group.

\M

\N {\bf Acknowledgements.} The first author gratefully acknowledges
support from Northwestern University for a three month visit in 2002.
The second author gratefully acknowledges partial support from l'
Universit\'e Louis Pasteur during a visit in academic year 1998-99
during which this collaboration began.

\BB

\section{The Braid group and Steinberg group associated to a root
system} \label{explm}

In discussing Coxeter groups and root systems, we will adopt the
terminology of Humphreys \cite[\S \S 1.2; 2.2; 2.8; 2.9]{H}, in which
all root systems are \emph{reduced} --- that is, the only root
parallel to $\al$ is $-\al$ --- but are not necessarily
crystallographic (\cite[2.9]{H}). Let $\Phi$ be an irreducible root
system and let $\DD$ be a fixed simple subsystem of $\Phi$. For $\al,
\be\in \DD$, define the integer $m(\al,\be)$ by

$$m(\al,\be) = \begin{cases}  1 & \text{ if }\alpha = \beta, \\
2 & \text{ if }\al \text{ and } \be \text{ are not connected in the
Coxeter graph  of } \Phi\\
m & \text{ if the edge connecting }\al \text{ and }
\be \text{ in the
Coxeter graph is labelled } m > 2
\end{cases}$$
(\emph{cf.} \cite{C}, \cite{H}, \cite{Bbki}).

The Weyl group, $W(\Phi)$, associated to $\Phi$ is presented with
generators the simple reflections $\{\ss_{\al} \mid \al\in \DD\}$ and
defining relations

\begin{equation} \label{defnweyl}
(\ss_{\al}\ss_{\be})^{m(\al,\be)}=1.
\end{equation}

If the root system $\Phi$ is \emph{crystallographic} and irreducible
 (\emph{i.e.} the root system of a simple Lie algebra over the complex
 numbers), we must have $m(\al,\be) = 1,2,3,4, \text{ or } 6$
 (\cite[Proposition, p. 38]{H}). In this case, there can be at most 2
 root lengths, and $W(\Phi)$ acts transitively on roots of the same
 length.

\subsection{The braid group of a root system} \label{brpresent}

For $\DD \subset \Phi$ (not necessarily crystallographic) as above,
define $Br(\Phi)$, the {\it braid group of type $\Phi$}, as the group
with generators $\{y_{\al} \mid \al\in \DD\}$ and defining relations
\begin{equation*}
(1.1.1)_{(\al,\be)} \qquad y_{\al} y_{\be} y_{\al}\dots = y_{\be} y_{\al}
       y_{\be}\dots
\end{equation*}

\noindent where there are $m=m(\al,\be)$ factors on each side.  The
map $y_\alpha \mapsto \sigma_\alpha$ defined on generators induces a
surjective homomorphism $Br(\Phi) \rightarrow W(\Phi)$ whose kernel,
$PBr(\Phi)$, is called the \emph{pure braid group of type }$\Phi$. We
will write $\bar{b}$ to denote the image in $W(\Phi)$ of an element $b
\in Br(\Phi)$.  Since $W(\Phi)$ is obtained from $Br(\Phi)$ by adding
the relations $y_\alpha^2 = 1 \text{ for } \alpha \in \Delta$, it
follows that $PBr(\Phi)$ is generated by all $y_\alpha^2, \alpha \in
\Delta$, together with their conjugates by arbitrary elements of
$Br(\Phi)$.

\subsection{The parametrized braid group of a root system} \label{explrelns}

Now assume that the root system $\Phi$ is crystallographic (so that
$m(\al,\be) = 2,3,4, \text{ or } 6$ when $\alpha \neq \beta$), and let
$A$ be a ring, assumed commutative if $\Phi \neq A_n$.
The {\it parametrized braid
group} of type
$\Phi$ with parameters in $A$, $Br(\Phi, A)$, has generators
$\y{\al}{a}, \al\in \DD$, $a\in A$, and relations

\begin{align*}
(A1)&&\y{\al}{a} \y{\al}{0} \y{\al}{b} &= \y{\al}{0} \y{\al}{0}
\y{\al}{a+b}  \\
(A1\times A1)&& \y{\al}{a} \y{\be}{b}&=
\y{\be}{b} \y{\al}{a}&\text{ if } m(\al,\be)=2\\
(A2) &&
\y{\al}{a}\y{\be}{b}\y{\al}{c} &= \y{\be}{c}\y{\al}{b+ac}\, \y{\be}{a}
&\text{ if } m(\al,\be)=3\\
(B2) &&y_{\alpha}^{a}y_{\beta}^{b}y_{\alpha}^{c}y_{\beta}^{d} &=
y_{\beta}^{-d}y_{\alpha}^{c-ad}y_{\beta}^{-b-2ac+a^2d}y_{\alpha}^{a}
&\text{ if }
m(\al,\be)=4
%\\
%(G2)&& \y{\al}{a}\y{\be}{b}\y{\al}{c}\y{\be}{d}\y{\al}{e}\y{\be}{f}&
%=\y{\be}{f}\y{\al}{e+af}\y{\be}{d-fb-3a^2e-3ea^2f}\y{\al}{c+2ae+a^2f}
%y{\be}{b+3ac+3a^2e+a^3f}\y{\al}{a}  \\ &&&&\text{ if } m(\al,\be)=6
%text{ and } \alpha < \beta
\end{align*}

\noindent
for $\al, \be\in \DD,\, \al \leq \beta$ (see below), and
$a,b,c,d,e,f\in A$.  There is also a relation $(G2)$ when
$m(\al,\be)=6$ which has 6 terms on each side; that relation and its
consequences will be discussed in \cite{LS2}.

The conditions $\alpha \leq \beta$ in relations $(A2) \text{ and }
(B2)$ refer to an ordering which will be made explicit in those cases
(\emph{cf.} \ref{brdef} and \ref{brDdef}).

When $\Phi = A_{n-1}$ this is precisely the group $Br_n(A)$ defined in
the introduction. The relations have been named to reflect the type of
irreducible crystallographic root system in which they occur.

Of course the root system $\Phi$ determines which of these relations
occur; for instance for the simply-laced root systems $\Phi=A_n, D_n,
E_6, E_7, E_8$ only the relations $(A1),(A1 \times A1)$ and $(A2) $
occur. In the cases $\Phi = B_n, C_n, F_4$ one needs to add relation
$(B2)$. The case $\Phi= G_2$ uses relations $(A1)$ and $(G2)$.

This presentation depends on the choice of a simple subsystem $\DD
\subset \Phi$, as do the presentations given in \S
\ref{brpresent}. However, it is easily checked that changing the
simple subsystem $\Delta$ results in a group isomorphic to the
original group.

When $A=0$, relation $(A1)$ becomes a tautology and the presentations
of $Br(\Phi,0)$ and $Br(\Phi)$ (\S \ref{brpresent}) coincide.
The unique ring homomorphism $A \rightarrow \{0\}$ induces a split
epimorphism $\pi : Br(\Phi, A) \rightarrow Br(\Phi)$, and thus $Br(\Phi,A)$
is isomorphic to the semi-direct product $$\Ker \pi  \rtimes Br(\Phi)\ .$$
Our aim is to identify $\Ker \pi $.

\subsection{The Steinberg group of a root system} \label{stbg}

Let $\Phi$ be an irreducible crystallographic root
system, and let $A$ be a commutative ring.  The {\it Steinberg group},
$St(\Phi,
A)$, of type $\Phi$ over $A$ (\emph{cf.} \cite{Sb} \cite{St}), is presented
with
generators $\{ x_{\al}(a)= x_{\al}^a \mid \al \in \Phi , a \in A\}$ (we
will freely
mix the notations $x_{\al}(a)$ and $ x_{\al}^a$ depending on how
complicated the coefficient $a$ is) subject to the relations

\begin{align*}& (R1)_\al&x_{\al}(a) x_{\al}(b) =&
x_{\al}(a+b)&\\& (R2)_{\al,\be}&[ x_{\al}(a) , x_{\be}(b)] =& \prod
x_{i\al+j\be}(N_{\al,\be,i,j}a^ib^j)&\\
\end{align*}
for all $\al,\be\in \Phi$ such that $\al+\be\neq 0$, where the product
is taken over all roots of the form $i\al+j\be, i,j\in
\mathbb{Z}_{>0}$ in some fixed order, and the $N_{\al,\be,i,j}$ are
certain integers depending only on the structure constants of the
simple Lie algebra of type $\Phi$.

\M We shall see that the Weyl group $W(\Phi)$ (and therefore, the
braid group $Br(\Phi)$) acts on $St(\Phi, A)$.  This action takes the
form
\begin{align*}  \label{action}
&&\sigma_\delta(x_\gamma(t)) &= x_{\sigma_\delta(\gamma)}(\eta(\delta,
\gamma)t)& \delta \in \Delta, \gamma \in \Phi & \end{align*} \noindent
where the $\eta(\delta, \gamma) = \pm 1$ are signs to be
determined. Thus we can form the semi-direct product $St(\Phi, A)
\rtimes Br(\Phi)$. The explicit form of relations \begin{math}(R1)_\al
\text{ and } (R2)_{\al,\be} \end{math} and of the $N_{\al,\be,i,j}$
(\cite[Theorem 5.2.2, p. 77; \S 4.2, p. 55]{C}) play an important
role in defining the action of the Weyl group on $St(\Phi, A)$.

In many cases, the group $St(\Phi, A)$ is the universal central
extension of $E(\Phi, A)$, the elementary subgroup of the points in
$A$ of a Chevalley-Demazure group scheme with (crystallographic) root
system $\Phi$ (\cite{St}).

\subsection{The Main Theorem and How It Is Proved} \label{mainthmpf}

In this paper and its sequel \cite{LS2}, we shall prove the
following theorem.

\begin{theorem}\label{mainthm}
For any crystallographic root system $\Phi$ and any commutative ring
$A$, the map
$$\phi : Br(\Phi,A) \to  St(\Phi, A) \rtimes  Br(\Phi)$$
induced by $\y{\al}{a} \mapsto (\x{\al}{a}, \y{\al}{})$ is an
isomorphism of groups.  This is true when $\Phi = A_{n-1}$ even when
$A$ is neither commutative nor unital.
\end{theorem}

Because the details in the proof of Theorem \ref{mainthm}
can become rather complicated, we summarize the steps first for the
reader's convenience.

\M \N $\bullet$  {\it Step (a).} We show that there exists a well-defined
map  $\phi : Br(\Phi,A)
\to St(\Phi, A) \rtimes Br(\Phi)$ induced by $\y{\al}{a} \mapsto (\x{\al}{a},
\y{\al}{})$ and we verify that  it is
a group homomorphism.

\M \N $\bullet$ {\it Step (b).} A crucial technical result needed to
show that $\phi$ is invertible is

\begin{lemma}[Pure
Braid Lemma] \label{PBL} For any $\alpha \in \Delta$ and any $a \in
A$, the element $y_\alpha^a(y_\alpha^0)^{-1} \in Br(\Phi, A)$ commutes
with every element of the pure braid group $PBr(\Phi)$ (viewed as a
subgroup of $Br(\Phi, A)$).
\end{lemma}

The proof of this lemma requires a manageable set of generators of the
pure braid group $PBr(\Phi)$. Such a result is classical for $\Phi=
A_{n-1}$ (\emph{cf.} \cite{Birman}), but less well-known in the other
cases. For the case $D_n$ in this paper, we deduce a set of generators
from the work of Digne and Gomi \cite{DG}.

\M \N $\bullet$ {\it Step (c).} We want to construct a group
homomorphism $\psi : St(\Phi, A) \rtimes Br(\Phi)\to Br(\Phi, A)$
which is inverse to $\phi$. We begin by defining a map $\psi$ which,
for $\alpha\in \DD$, sends $y_\alpha$ to $y_\alpha ^0$ and $x_\alpha
^a$ to $y_\alpha^a(y_\alpha^0)^{-1}$. But $St(\Phi, A)$ has many more
generators, namely all the elements $x_\alpha ^a$ for $\alpha \in
\Phi$ (and not just those in the simple subsystem $\DD$), so we need
to extend the definition of $\psi(x_\alpha ^a)$ to any $\alpha \in
\Phi$. Now since the Weyl group acts transitively on the set of roots
of each length, we can choose $\ss\in W(\Phi)$ such that
$\ss(\alpha)\in \DD$. Choosing a lifting $\tilde \ss \in Br(\Phi)$ of
$\ss$, we define $\psi( x_\alpha ^a)$ as $\tilde \ss\mm \cdot
\psi(x_{\ss(\alpha)}^a)$. To prove that $\psi$ well-defined -- that
is, independent of the choice of $\sigma$ and of its lifting $\tilde
\sigma$ -- and is a group homomorphism requires \emph{Step (b)}.

Since it is clear that $\phi\circ \psi$ and $\psi\circ \phi $ are
identity maps, these three steps establish the main theorem.  \M

\M In this paper we treat completely two cases for which $m(\alpha,
\beta) \leq 3$, those of $\Phi=A_{n-1}$ (which does not require the
hypothesis of commutativity on $A$) and $\Phi=D_n$.  

\M Before beginning the case by case proof we state two results which
will prove useful later on.  \M

\begin{lemma}\label{exrelns} The following relations in $Br(\Phi,
A)$ are consequences of relation $(A1)$ of \S  \ref{explrelns}.
$$\begin{array}{rcl}
(\y{\alpha}{0} \y{\alpha}{0} )\y{\alpha}{a} &=& \y{\alpha}{a} (\y{\alpha}{0} \y{\alpha}{0} ) ,\\
\y{\alpha}{a} \ym{\alpha}{0} \y{\alpha}{b}   &=& \y{\alpha}{a+b} ,\\
\y{\alpha}{-a}              &=& \y{\alpha}{0} \ym{\alpha}{a} \y{\alpha}{0} .\\
\end{array}$$
\end{lemma}

\N {\it  Proof.} Replacing  $b$  by 0  in  $(A1)$ shows  that
$\y{\alpha}{0} \y{\alpha}{0}$ commutes with $\y{\alpha}{a}$. From this
follows
$$\begin{array}{rcl} \y{\alpha}{a} \ym{\alpha}{0} \y{\alpha}{b} &=&
(\y{\alpha}{0} \y{\alpha}{0})^{-1} \y{\alpha}{a} (\y{\alpha}{0}
\y{\alpha}{0}) (\y{\alpha}{0})^{-1} \y{\alpha}{b} \\ &=&
(\y{\alpha}{0} \y{\alpha}{0})^{-1} \y{\alpha}{a} \y{\alpha}{0}
\y{\alpha}{b} \\ &=& (\y{\alpha}{0} \y{\alpha}{0})^{-1} \y{\alpha}{0}
\y{\alpha}{0} \y{\alpha}{a+b}\\ &=& \y{\alpha}{a+b} . \\
\end{array}$$
Putting $b=-a$ in the second relation yields the third relation.\hfill
\qed \M

\begin{lemma}\label{conjugation} Assume that the Pure Braid Lemma
\ref{PBL} holds for the root system $\Phi$. Suppose that $\alpha \in
\Delta$ and $b\in Br(\Phi)$ are such that $\bar{b} (\alpha)\in \Delta$
(where $\bar{b}\in W(\Phi)$ denotes the image of $b$; \emph{cf.}
\S \ref{brpresent}). Then for any $a\in A$
$$b \y{\alpha}{a}(\y{\alpha}{0})^{-1} b^{-1} = \y{\bar{b}
(\alpha)}{a}(\y{\bar{b} (\alpha)}{0})^{-1}\quad   \hbox {in }
Br(\Phi, A).$$
\end{lemma}

\N {\it Proof.} First let us show that that there exists ${b'}\in
Br(\Phi)$ such that $b' \y{\alpha}{a} {b'}^{-1} = \y{\bar{b}
(\alpha)}{a}$. The two roots $\alpha$ and $\bar{b} (\alpha)$ have the
same length, hence they are connected, in the Dynkin diagram, by a
finite sequence of edges with $m=3$ (\cite[Lemma 3.6.3]{C}). Therefore
it is sufficient to prove the existence of ${b'}$ when $\alpha$ and
$\bar{b} (\alpha)$ are adjacent.  In that case $\alpha$ and $\bar{b}
(\alpha)$ generate a subsystem of type $A_2$; we may assume $\alpha =
\alpha_1$ and $\bar{b} (\alpha) = \alpha_2 \in A_2$; and we can use
the particular case of relation $(A2)$, namely
$$\y{\alpha_1}{0}\y{\alpha_2}{0}\y{\alpha_1}{a}=\y{\alpha_2}{a}\y{\alpha_1}{0}\y{\alpha_2}{0}$$
to show that
$$\y{\alpha_1}{0}\y{\alpha_2}{0}\y{\alpha_1}{a}(\y{\alpha_2}{0})^{-1}(\y{\alpha_1}{0})^{-1}=\y{\alpha_2}{a}.$$
(Here $b' = \y{\alpha_1}{0}\y{\alpha_2}{0}$, $\alpha = \alpha_1 =
-\epsilon_1 + \epsilon_2, \bar{b'} (\alpha)= \alpha_2 =
-\epsilon_2+\epsilon_3$.)

To conclude the proof of the Lemma it is sufficient to show that
$$b \y{\alpha}{a}(\y{\alpha}{0})^{-1} b^{-1} =
\y{\alpha}{a}(\y{\alpha}{0})^{-1}$$ whenever $b(\alpha) =
\alpha$. According to \cite[Theorem, p. 22]{H}, $\bar{b}$ is a product
of simple reflections $\sigma_{\alpha_i}$ for $\alpha_i \in \Delta$
which are not connected to $\al$ in the Dynkin diagram of
$\Delta$. Hence we can write $b$ as the product of an element in the
pure braid group and generators $y_{\al_i}$ which commute with
$\y{\alpha}{a}$ by relation $(A1\times A1)$.  Since we have assumed
that the Pure Braid Lemma holds for $\Phi$, we can thus conclude that
$b \y{\alpha}{a}(\y{\alpha}{0})^{-1} b^{-1} =
\y{\alpha}{a}(\y{\alpha}{0})^{-1}$ as desired.  \hfill \qed \M

\section{The parametrized braid group  for $\Phi = A_{n-1}$}
\label{explpix}

Let $A$ be a ring, in this section not necessarily unital or
commutative. We consider the special case of the parametrized braid
group $Br(\Phi, A)$ for $\Phi = A_{n-1}$. In this case we will write
$y_i$ instead of $y_{\alpha_i}$, and $Br_n(A)$ for $Br(A_{n-1},
A)$. Because $A_{n-1}$ is a simply-laced root system, we have

\begin{definition}\label{brdef} Let $A$ be a ring (not necessarily unital nor
commutative). The {\it parametrized braid group} $Br_n(A)$ is
generated by the elements $\y{i}{a}$, where $i$ is an integer, $1\leq
i\leq n-1$, and $a$ is an element of $A$, subject to the relations

\begin{align*}
(A1)\quad&& \y{i}{a} \y{i}{0} \y{i}{b}&= \y{i}{0} \y{i}{0} \y{i}{a+b}
& \\
(A1 \times A1)\quad&& \y{i}{a} \y{j}{b} &= \y{j}{b}
\y{i}{a}&\text{ if } |i-j|\geq 2,\\
(A2)\quad&&
\y{i}{a}\y{i+1}{b}\y{i}{c}&= \y{i+1}{c}\y{i}{b+ac}\, \y{i+1}{a} & \\
\end{align*}
for any $a,b,c \in A$.
\end{definition}

The geometric motivation for the defining relations of this group, and
its connection with braids, can be seen in the following figures in
which $u, v,w$ (the colors) are elements of $A$, and $a,b,c$ are the
{\it coefficients of spilling}. Relation $(A1)$ comes from Figure 2.

%figure2
\setlength{\unitlength}{0.8cm}

\centerline {
\begin{picture}(20,16)
\put(0.9,15.5){$u$}
\put(2.9,15.5){$v$}
\put(1,15){\line(0,-1){1}}
\put(3,15){\line(0,-1){1}}
\put(1,12){\line(0,-1){1}}
\put(3,12){\line(0,-1){1}}
\put(1,14){\line(1,-1){2}}
\put(3,14){\line(-1,-1){0.8}}
\put(1,12){\line(1,1){0.8}}
\put(0,12.7){$y_i^b$}
\put(2.9,10.5){$u$}
\put(0.3,10.5){$v+bu$}
\put(1,10){\line(0,-1){1}}
\put(3,10){\line(0,-1){1}}
\put(1,7){\line(0,-1){1}}
\put(3,7){\line(0,-1){1}}
\put(1,9){\line(1,-1){2}}
\put(3,9){\line(-1,-1){0.8}}
\put(1,7){\line(1,1){0.8}}
\put(0,7.7){$y_i^0$}
\put(0.9,5.4){$u$}
\put(2.2,5.4){$v+bu$}
\put(1,5){\line(0,-1){1}}
\put(3,5){\line(0,-1){1}}
\put(1,2){\line(0,-1){1}}
\put(3,2){\line(0,-1){1}}
\put(1,4){\line(1,-1){2}}
\put(3,4){\line(-1,-1){0.8}}
\put(1,2){\line(1,1){0.8}}
\put(0,2.7){$y_i^a$}
\put(2.9,0.5){$u$}
\put(-1,0.5){$v+bu+au$}
\put(6,7.7){$=$}
%partie droite de la figure
\put(8.9,15.5){$u$}
\put(10.9,15.5){$v$}
\put(9,15){\line(0,-1){1}}
\put(11,15){\line(0,-1){1}}
\put(9,12){\line(0,-1){1}}
\put(11,12){\line(0,-1){1}}
\put(9,14){\line(1,-1){2}}
\put(11,14){\line(-1,-1){0.8}}
\put(9,12){\line(1,1){0.8}}
\put(13,12.7){$y_i^{a+b}$}
\put(7,10.5){$v+(a+b)u$}
\put(10.9,10.5){$u$}
\put(9,10){\line(0,-1){1}}
\put(11,10){\line(0,-1){1}}
\put(9,7){\line(0,-1){1}}
\put(11,7){\line(0,-1){1}}
\put(9,9){\line(1,-1){2}}
\put(11,9){\line(-1,-1){0.8}}
\put(9,7){\line(1,1){0.8}}
\put(13,7.7){$y_i^0$}
\put(8.9,5.4){$u$}
\put(10,5.4){$v+(a+b)u$}
\put(9,5){\line(0,-1){1}}
\put(11,5){\line(0,-1){1}}
\put(9,2){\line(0,-1){1}}
\put(11,2){\line(0,-1){1}}
\put(9,4){\line(1,-1){2}}
\put(11,4){\line(-1,-1){0.8}}
\put(9,2){\line(1,1){0.8}}
\put(13,2.7){$y_i^0$}
\put(7,0.5){$v+(a+b)u$}
\put(10.9,0.5){$u$}
\end{picture}
}

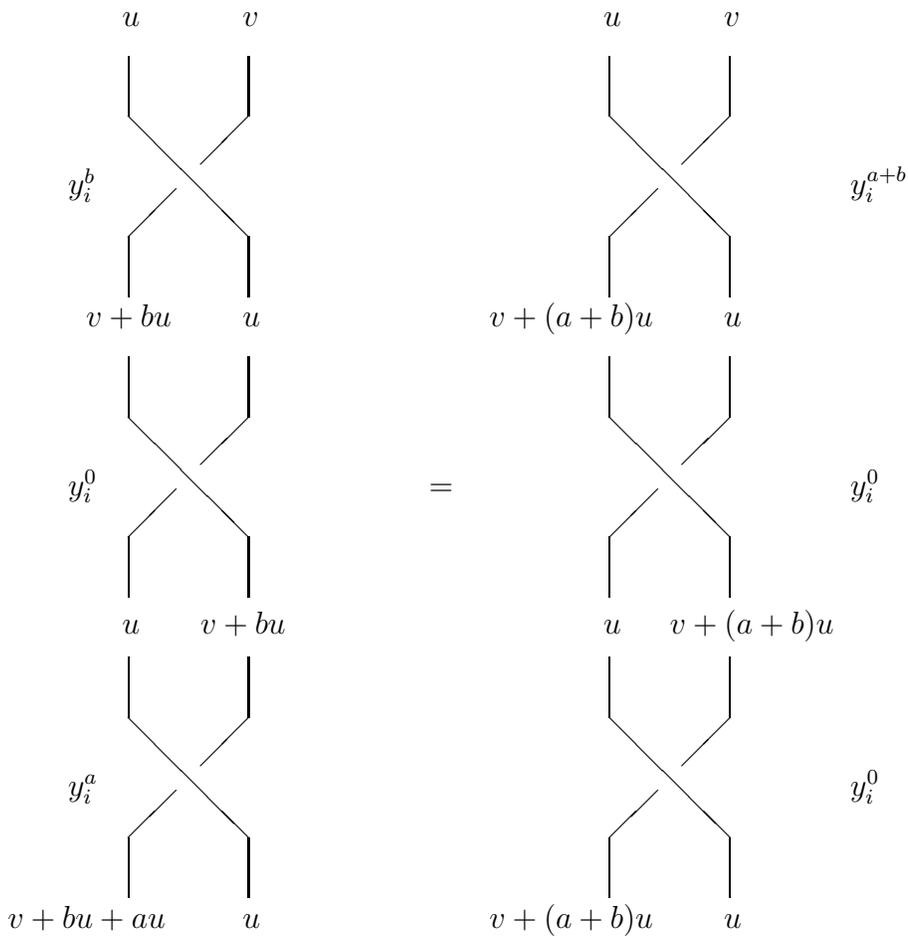
\centerline {Figure 2}
\B

Relation $(A1\times A1)$ arises because the actions of $\y{i}{a}$ and
of $\y{j}{b}$ on the strands of the braid are disjoint when $\vert
i-j\vert \geq 2$, so that these two elements commute.

Relation $(A2)$ derives from Figure 3.

%figure3
\centerline{
\begin{picture}(20,16)
\put(0.9,15.5){$u$}
\put(2.9,15.5){$v$}
\put(4.9,15.5){$w$}
\put(1,15){\line(0,-1){1}}
\put(3,15){\line(0,-1){1}}
\put(1,12){\line(0,-1){1}}
\put(3,12){\line(0,-1){1}}
\put(1,14){\line(1,-1){2}}
\put(3,14){\line(-1,-1){0.8}}
\put(1,12){\line(1,1){0.8}}
\put(5,15){\line(0,-1){4}}
\put(0,12.7){$y_i^c$}
\put(0.5,10.5){$v+cu$}
\put(2.9,10.5){$u$}
\put(4.9,10.5){$w$}
\put(3,10){\line(0,-1){1}}
\put(5,10){\line(0,-1){1}}
\put(3,7){\line(0,-1){1}}
\put(5,7){\line(0,-1){1}}
\put(3,9){\line(1,-1){2}}
\put(5,9){\line(-1,-1){0.8}}
\put(3,7){\line(1,1){0.8}}
\put(1,10){\line(0,-1){4}}
\put(0,7.7){$y_{i+1}^b$}
\put(0.5,5.6){$v+cu$}
\put(2.2,5.6){$w+bu$}
\put(4.9,5.6){$u$}
\put(1,5){\line(0,-1){1}}
\put(3,5){\line(0,-1){1}}
\put(1,2){\line(0,-1){1}}
\put(3,2){\line(0,-1){1}}
\put(1,4){\line(1,-1){2}}
\put(3,4){\line(-1,-1){0.8}}
\put(1,2){\line(1,1){0.8}}
\put(5,5){\line(0,-1){4}}
\put(0,2.7){$y_i^a$}
\put(0,0.5){$w+bu$}
\put(0,0){$+a(v+cu)$}
\put(2.5,0.5){$v+cu$}
\put(4.9,0.5){$u$}
\put(6.5,7.7){$=$}
%partie droite de la figure
%
\put(7.9,15.5){$u$}
\put(9.9,15.5){$v$}
\put(11.9,15.5){$w$}
\put(10,15){\line(0,-1){1}}
\put(12,15){\line(0,-1){1}}
\put(10,12){\line(0,-1){1}}
\put(12,12){\line(0,-1){1}}
\put(10,14){\line(1,-1){2}}
\put(12,14){\line(-1,-1){0.8}}
\put(10,12){\line(1,1){0.8}}
\put(8,15){\line(0,-1){4}}
\put(13,12.7){$y_{i+1}^{a}$}
\put(7.9,10.5){$u$}
\put(9.5,10.5){$w+av$}
\put(11.9,10.5){$v$}
\put(8,10){\line(0,-1){1}}
\put(10,10){\line(0,-1){1}}
\put(8,7){\line(0,-1){1}}
\put(10,7){\line(0,-1){1}}
\put(8,9){\line(1,-1){2}}
\put(10,9){\line(-1,-1){0.8}}
\put(8,7){\line(1,1){0.8}}
\put(12,10){\line(0,-1){4}}
\put(13,7.7){$y_i^{b+ac}$}
\put(7.4,5.5){$w+av$}
\put(7.4,5){$+(b+ac)u$}
\put(9.9,5.5){$u$}
\put(11.9,5.6){$v$}
\put(10,5){\line(0,-1){1}}
\put(12,5){\line(0,-1){1}}
\put(10,2){\line(0,-1){1}}
\put(12,2){\line(0,-1){1}}
\put(10,4){\line(1,-1){2}}
\put(12,4){\line(-1,-1){0.8}}
\put(10,2){\line(1,1){0.8}}
\put(8,4.8){\line(0,-1){3.8}}
\put(13,2.7){$y_{i+1}^c$}
\put(7.4,0.5){$w+av$}
\put(7.4,0){$+(b+ac)u$}
\put(9.4,0.5){$v+cu$}
\put(11.9,0.5){$u$}
\end{picture}
}

\smallskip
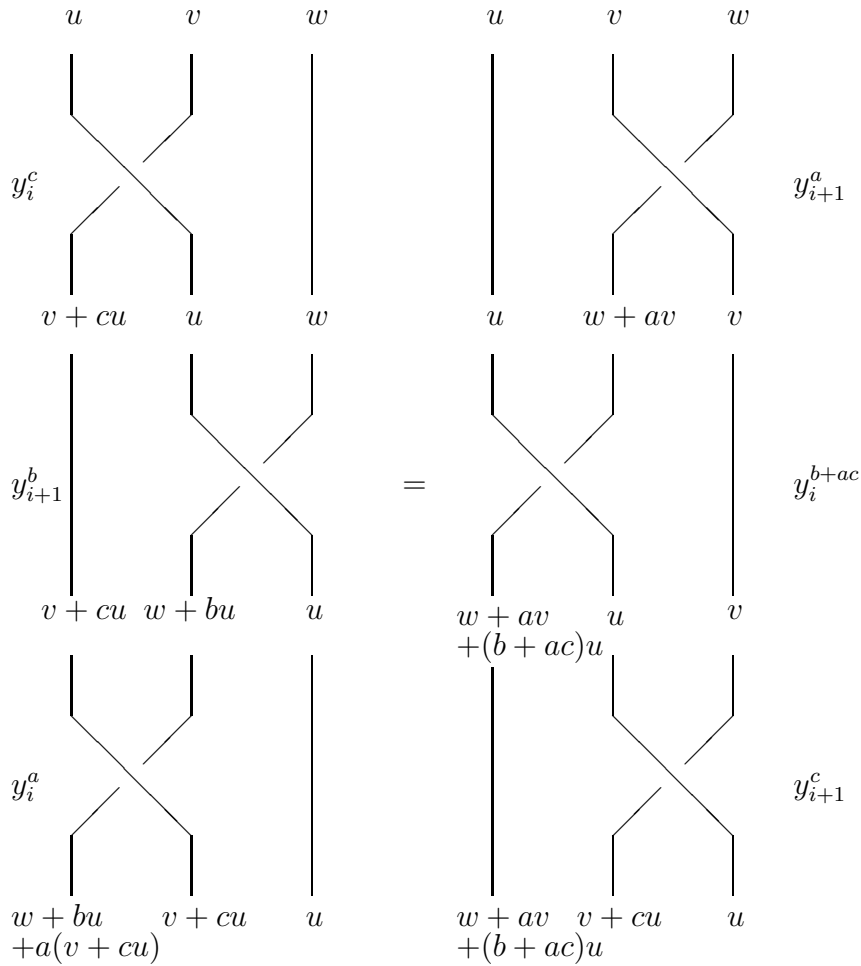
\centerline {Figure 3}
\B

\subsection{Braid group and pure braid group} \label{simplnot}

  The group
$Br_n (0) = Br_n$ is the classical Artin braid group with generators
$y_i\, , 1\leq i\leq n-1,
$ and relations
$$\begin{array}{rcl}
y_i\,  y_j               &=& y_j\, y_i\,  , \qquad |i-j| \geq 2,\\
y_i\,  y_{i+1}\,  y_i\,  &=&  y_{i+1}\,  y_i\,  y_{i+1}\, .\\
\end{array}$$

The quotient of $Br_n$ by the relations $y_i y_i = 1 \, , 1 \leq i
\leq n-1$ is the symmetric group ${\mathcal{S}_n}$; the image of $b
\in Br_n$ in ${\mathcal{S}_n}$ is denoted by $\bar{b}$. The kernel of
the surjective homomorphism $Br_n \to {\mathcal{S}_n}$ is the {\it
pure braid group}, denoted $PBr_n$. It is generated by the elements
$$ {\aa}_{j,i}:= y_j y_{j-1} \cdots y_i y_i \cdots y_{j-1}y_j ,$$
 for $n\geq j\geq i\geq 1$, (\cite{Birman}; see Figure 4 below).

%purebraid
\begin{picture}(16,8)
\put(1,0){\line(0,1){7}}
\put(2,0){\line(0,1){2.8}}
\put(2,3.2){\line(0,1){3.8}}
\put(3,0){\line(0,1){2.3}}
\put(3,2.7){\line(0,1){4.3}}
\put(5,0){\line(0,1){1.3}}
\put(5,1.7){\line(0,1){5.3}}
\put(6,0){\line(0,1){0.8}}
\put(6,1.2){\line(0,1){5.8}}
\put(8,0){\line(0,1){7}}

\qbezier(7,0)(7,0.5)(6,1)
\put(6,1){\line(-2,1){4}}
\qbezier(2,3)(1,3.5)(1.8,3.9)
\put(2.2,4.1){\line(2,1){0.6}}
\put(3.2,4.6){\line(2,1){0.6}}
\put(5.2,5.6){\line(2,1){0.6}}
\qbezier(6.2,6.1)(7,6.5)(7,7)

\put(1.9,7.3){$i$}
\put(5.9,7.3){$j$}
\put(3.5,3.4){$\cdots$}
\end{picture}

\smallskip
\centerline{Figure 4: the pure braid $ {\aa}_{j,i}$}
\B

\begin{lemma}[Pure Braid Lemma for $A_{n-1}$] \label{pblan} Let $y_k^a$
be a generator of $Br_n(A)$ and let $\oo \in PBr_n=
PBr_n(0)$. Then there exists $\oo '\in Br_n$, independent of $a$,  such
that
$$ y_k^a\ \oo = \oo '\ y_k^a.$$ Hence for any integer $k$ and any
element $a\in A$, the element $y_k^a(y_k^0)\mm\in PBr_n(A)$ commutes
with every element of the pure braid group $PBr_n$. \end{lemma}

\M
\N {\bf Notation.} Before beginning the proof of Lemma \ref{pblan}, we want to
simplify our notation.  We will abbreviate

$$\begin{array}{rcl}
y_k^0 &=& \kk \\

(y_k^0)^{-1}  &=& \kk^{-1} \\

y_k^a &=& \kk^a
\end{array}$$
Note that $\kk^{-1}$ does {\it not} mean $\kk^a$ for $a=-1$.

If there exist $\omega '$ and $\omega ''\in Br_n$ such that $\kk^a
\omega = \omega ' \jj^a \omega ''$, we will write
$$\begin{array}{rcl}\kk^a \omega \sim \jj^a \omega '' .\end{array}$$
Observe that $\sim$ is {\it not} an equivalence relation, but it is
compatible with multiplication on the right by elements of $Br_n$: if
$\eta \in Br_n$,
$$\kk^a \omega \sim \jj^a \omega '' \Leftrightarrow \kk^a\omega = \omega '
\jj^a \omega '' \Leftrightarrow \kk^a\omega \eta = \omega ' \jj^a \omega ''
\eta
 \Leftrightarrow \kk^a \omega\eta \sim \jj^a \omega '' \eta$$

\N For instance
$$\begin{array}{rclr}
\kk^a\ \kk\ \kk &\sim & \kk^a & \hbox{by Lemma \ref{exrelns}},\\
\kk^a \ \kmkm \ \kk & \sim & \kmkm ^a &  \hbox{by relation (A2)},\\
\kk^a \ \kmkm \  & \sim & \kmkm ^a\ \kk\ &   \hbox{by relation (A2)}.\\
\end{array}$$

\M \N {\it Proof of Lemma \ref{pblan}.}  It is clear that the first
assertion implies the second one: since $\oo'$ is independent of $a$,
we can set $a=0$ to determine that $\oo' = y_k^0\ \oo \
(y_k^0)^{-1}$. Substituting this value of $\oo'$ in the expression $
y_k^a\ \oo = \oo '\ y_k^a$ completes the proof.

\M To prove the first assertion, we must show, in the notation just
introduced, that $\kk^a \omega \sim \kk^a$ for every $\omega \in
PBr_n$, and it suffices to show this when $\omega$ is one of the
generators ${\aa}_{j,i}= \jj\ \jmjm\ \cdots \ii\ \ii \cdots \ \jmjm\
\jj$ above. Since $y_k^a$ commutes with $y_i^0$ for all $i\neq k-1, k,
k+1$, it commutes with ${\aa}_{j,i}$ whenever $k< i-1$ or whenever $k
> j+1$. So we are left with the following 3 cases:
\begin{subequations}
\begin{eqnarray}
\qquad \quad \jpjp^a {\aa}_{j,i} &\sim& \jpjp^a  \label{relna1} \\
\qquad \quad \jj^a {\aa}_{j,i} &\sim& \jj^a  \label{relna2} \\
\qquad \quad \kk^a {\aa}_{j,i} &\sim& \kk^a  \quad i-1 \leq k \leq j-1
\label{relna3}
\end{eqnarray}
\end{subequations}

\noindent Our proof of case \eqref{relna1} is by induction on the
half-length of $\oo = {\aa}_{j,i}$. When $i=k-1 $, $\oo = \kmkm
\ \kmkm $, and we have
$$\begin{array}{rcl}
\kk^a \ \kmkm  \ \kmkm  & =   &\kk^a \ \kmkm \ \kk\mm \kk\ \kmkm \\
                &\sim &\kmkm ^a\ \kk \ \kmkm \\
                &\sim &\kk^a\ .\\
\end{array}$$
and more generally,
$$\begin{array}{rcl}
\kk^a \ \kmkm \ \kmmkmm \cdots \kmmkmm\ \kmkm &\sim& \kmkm ^a\ \kk \
\kmmkmm  \cdots \kmmkmm \kmkm \\
&=& \kmkm ^a \ \kmmkmm  \cdots  \kmmkmm  \ \kk\ \kmkm  \\
&\sim& \kmkm ^a \ \kk\  \kmkm \qquad \hbox{(by induction)}, \\
&\sim& \kk^a\ .\\
\end{array}$$

\N The proof of case \eqref{relna2} is also by induction on the
half-length of $\oo = {\aa}_{j,i}$. When $i=k$, $\oo = \kk\
\kk$, and we have
$$\kk^a \kk\ \kk \sim \kk^a \qquad \hbox {by Lemma \ref{exrelns}}.$$
Then
$$\begin{array}{rcl}
 \kk^a \kk\ \kmkm  \cdots \kmkm  \ \kk &   = &\kk^a \kk\ \kmkm \ \kk\
\kk\mm \ \kmmkmm  \cdots \kmmkmm \  \kmkm  \ \kk  \\
  &\sim &\kk^a  \kmkm \ \kk\ \kmkm  \ \kmmkmm  \cdots \kmmkmm \ \kk\mm
\kmkm  \ \kk  \\
  &\sim &\kmkm ^a \kmkm \  \kmmkmm  \cdots \kmmkmm \ \kmkm  \ \kk \ \kmkm
\mm  \\
  &\sim &\kmkm ^a \ \kk \ \kmkm \mm \qquad \hbox{(by induction)} \\
  &\sim &\kk^a.\\
\end{array}$$

\N For case \eqref{relna3} it is sufficient to check the cases
\begin{subequations}
\begin{eqnarray}
\quad\oo & = & \kpkp \ \kpkp \label{relna4} \\
\quad \oo &=& \kpkp \ \kk\ \kk\ \kpkp \label{relna5} \\
\quad \oo &=& \kpkp \ \kk\ \kmkm  \cdots \kmkm  \ \kk\ \kpkp \label{relna6}
\end{eqnarray}
\end{subequations}

\N which are proved as follows:

$$\begin{array}{rcl} \kk^a \kpkp \ \kpkp & = & \kk^a \kpkp \ \kk\mm
\kk \ \kpkp\\ &\sim & \kpkp
^a \ \kk\ \kpkp \\ &\sim &\kk^a  \qquad \hbox { Case \eqref{relna4}}
\end{array}$$

$$\begin{array}{rcl} \kk^a \kpkp \kk\kk \kpkp & \sim &\kpkp^a \kk \
\kpkp\\ &\sim &\kk^a  \qquad \hbox { Case \eqref{relna5}}
\end{array}$$

$$\begin{array}{rcl} \kk^a \kpkp \ \kk\ \kmkm \cdots \kmkm \ \kk\
\kpkp & \sim &\kpkp^a \kmkm \cdots \kmkm \ \kk \ \kpkp \\ &\sim &\kpkp^a \kk\
\kpkp \\ &\sim
&\kk^a \qquad \hbox { Case
\eqref{relna6}}
\end{array}$$
\hfill \qed
\M

\begin{proposition} \label{conjan}
 For any $\omega \in Br_n(0)$ the element $ \omega y_k^a (y_k^0)^{-1}
 \omega\mm $ depends only on the class $\bar{\omega}$ of $\omega$ in
 ${\mathcal{S}_n}$. Moreover if $\bar{\omega}(j) = k$ and
$\bar{\omega}(j+1) = k+1
 $, then $ \omega y_k^a (y_k^0)\mm \omega\mm =y_j^a (y_j^0)\mm$.
\end{proposition}

\N {\it Proof.} The first statement is a consequence of Lemma
\ref{PBL} since the Weyl group $W(A_{n-1}) = {\mathcal{S}_n}$ is the
quotient of
$Br_n$ by $PBr_n$.

The second part is a consequence of relation $(A1\times A1)$ and
the following computation:
$$\begin{array}{rcl}
\22 \mm \11 \mm \22 ^a \11  \22  &=&\11 ^a \22 \mm \11 \mm \22 \mm \11  \22 \\
&=&\11 ^a \11 \mm \22 \mm \11 \mm \11  \22 \\
&=&\11 ^a \11 \mm\ .\\
\end{array}$$
{ }\hfill $\square$
\M

\subsection{The Steinberg group for $\Phi = A_{n-1}$ and an action of
  the Weyl group}

 When $\Phi = A_{n-1}$, the Steinberg group of \S \ref{stbg} is
well-known and customarily denoted $St_n(A)$. In that case it is
customary to write $\x{ii+1}{a} = x_{ii+1}(a)$ for the element $x_{\alpha}(a) , \alpha = \epsilon _i - \epsilon _{i+1} \in
\Delta$, and, more generally,
$\x{ij}{a}$ when $\al = \epsilon_i - \epsilon_j \in A_{n-1}$.

\begin{definition}\label{brAdef} The Steinberg group of the ring $A$,
denoted $St_n(A)$, is presented by
the generators
$\x{ij}{a} , 1\leq i,j\leq n, i\neq j, a \in A$ subject to the relations
\begin{align*}
(St0)\quad&&\x{ij}{a} \x{ij}{b} &= \x{ij}{a+b}                      &\\
(St1)\quad&&\x{ij}{a} \x{kl}{b} &= \x{kl}{b} \x{ij}{a} ,            &\quad
i\neq l, j\neq k \\
(St2)\quad&&\x{ij}{a} \x{jk}{b} &= \x{jk}{b} \x{ik}{ab} \x{ij}{a} , &\quad
i\neq k .\\
\end{align*}
\end{definition}
We should make two observations about this definition. First, it
follows from $(St0)$ that $\x{ij}{0} = 1$.  Second, relation $(St2)$
is given in a perhaps unfamiliar form (\emph{i.e.}  different from the
commutator relation $(R2)$ of \S \ref{stbg} in the case $\Phi =
A_{n-1}$). We have chosen this form, which is easily seen to be
equivalent to $(R2)$, because of its geometric significance
(\emph{cf}. \cite{KS} for the relationship with the Stasheff
polytope), and for the simplification it brings in computation. \M

The Weyl group $W(A_{n-1})$ is isomorphic to the symmetric group
$\mathcal{S}_n$. Its action on the Steinberg group is induced by the
formula
\begin{equation} \label{defnaction}
\sigma \cdot \x{ij}{a} := \x{ \sigma (i) \sigma (j)} a , \quad
\sigma \in \mathcal{S}_n , a \in A
\end{equation}
Hence, in the notation of \S 1.3, one has $\eta(\delta, \gamma) = +1$.

\subsection{The main result for $\Phi = A_{n-1}$}

\begin{theorem} \label{homom} For any (not necessarily unital) ring $A$ the
map
$$\phi : Br_n(A) \to St_n(A)\rtimes Br_n$$ from the parametrized braid
group to the semi-direct product of the Artin braid group with the
Steinberg group induced by $\phi (\y{i}{a}) = \x{i\, i+1}{a} \, y_i$
is a group isomorphism.
\end{theorem}

\N {\it Proof.}  {\it Step (a).} We show that $\phi$ is a well-defined
group homomorphism.  \M

\N $\bullet$ Relation (A1):
$$\begin{array}{rcl}
\phi (\y{i}{a} \y{i}{0} \y{i}{b} ) &= &\x{i\, i+1}{a} y_i\, y_i\,\x{i\,
i+1}{b} y_i\, ,\qquad {\rm since}\quad \x{i\, i+1}{0} = 1,\\
&= &y_i\,  y_i\,   \x{i\, i+1}{a} \x{i\, i+1}{b}  y_i\, \qquad {\rm
since}\quad \overline{y_i\,  y_i\, } = 1\in
{\mathcal{S}_n},\\ &= &y_i\,  y_i\,   \x{i\, i+1}{a+b}   y_i\,  \qquad
\qquad  {\rm by\  } (St0),\\
 &= &\phi (\y{i}{0} \y{i}{0} \y{i}{a+b} ).\\
\end{array}$$
 \N $\bullet$ Relation $(A1\times A1)$ follows immediately from $(St1)$.

\N $\bullet$ Relation $(A2)$ is proved by using the relations of $Br_n$ and
the 3 relations $(St0),(St1),(St2)$ as follows:
$$\begin{array}{rcl}
\phi(\y{i}{a}\y{i+1}{b} \y{i}{c}) &=&\x{i\, i+1}{a} y_i\,\x{i+1\,
i+2}{b}\uu{y_{i+1}\,\x{i\, i+1}{c}}y_i\,  \\
&=&\x{i \, i+1}{a} y_i\,  \uu{\x{i+1\, i+2}{b}  \x{i\, i+2}{c} } \,
y_{i+1}\,  y_i\, , \quad \\
&=&\x{i \, i+1}{a} \uu{y_i\,  \x{i \, i+2}{c} } \x{i+1\, i+2}{b}\,
y_{i+1}\,  y_i\, , \quad \\
&=&\uu{\x{i \, i+1}{a} \x{i+1 \, i+2}{c} }\uu{y_i\,   \x{i+1\, i+2}{b}} \,
y_{i+1}\,  y_i\, , \quad \\
&=&\x{i+1 \, i+2}{c}  \x{i \, i+2}{ac}  \uu{\x{i \, i+1}{a}  \x{i\, i+2}{b}
} \uu{y_i\,   y_{i+1}\,  y_i}\, , \quad \\
&=&\x{i+1 \, i+2}{c}  \x{i \, i+2}{ac+b}  \uu{\x{i \, i+1}{a}   y_{i+1}} \,
y_i\, y_{i+1}\,  , \quad \\
&=&\x{i+1 \, i+2}{c} \uu{ \x{i \, i+2} {ac+b} \, y_{i+1}}\, \uu{ \x{i \,
i+2}{a} \, y_{i}\, y_{i+1}}\, , \quad \\
&=&\uu{\x{i+1 \, i+2}{c} y_{i+1}}\,\uu{ \x{i \, i+1} {ac+b}
y_i}\,\uu{\x{i+1 \, i+2}{a} \, y_{i+1}}\, ,\quad \\
&=&\phi ( \y{i+1}{c} \, \y{i}{b+ac} \, \y{i+1}{a} ).\\
\end{array}$$

 {\it Step (b).} This is the Pure Braid Lemma \ref{pblan} for $A_n$.
\M

 {\it Step (c).} We construct a homomorphism $\psi: St_n(A)\rtimes
Br_n \to Br_n (A)$. We first construct $\psi: St_n(A)\to \Ker \pi$,
where $\pi$ is the surjection $Br_n(A) \to Br_n$, by setting
$$\begin{array}{rcl}
\psi ( \x{ij}{a} ) &:=& \oo \ \y{k}{a} \ym{k}{0} \ \oo ^{-1}\\
\end{array}$$
where $\oo$ is an element of $Br_n$ such that $\bar \oo (k) = i$ and
$\bar \oo (k+1) = j$ (for instance, $\psi (\x{12}{a})= \y{1}{a}
\ym{1}{0}$ and $\psi (\x{13}{a})= \y{2}{0} (\y{1}{a} \ym{1}{0})
\ym{2}{0}$). Observe that this definition does not depend on the
choice of $\bar\oo$ (by Lemma \ref{conjugation}), and does not depend
on how we choose a lifting $\oo$ of $\bar\oo$ (by the Pure Braid Lemma
\ref{pblan}).

In order to show that $\psi$ is a homomorphism, we must demonstrate
that the Steinberg relations are preserved.  \M

\N $\bullet $ Relation $(St0)$:  it suffices to show that  $\psi(\x{12}{a}
\x{12}{b} ) = \psi( \x{12}{a+b} )$,
$$\begin{array}{rlr}
\psi (\x{12}{a} \x{12}{b} ) &=\y{1}{a} \ym{1}{0}  \y{1}{b} \ym{1}{0} =
\y{1}{a} (\y{1}{0} )^{-2} \y{1}{0}
 \y{1}{b} \ym{1}{0}&\\
 &= (\y{1}{0} )^{-2}\y{1}{a} \y{1}{0} \y{1}{b} \ym{1}{0}& \hbox {
by \ref{pblan}},\\
&=\y{1}{a+b}  \ym{1}{0} &\hbox{by (A1)},\\
 &=\psi( \x{12}{a+b} ).& \\
\end{array}$$
\N $\bullet $ Relation $(St1)$: it suffices to show that
$\psi(\x{12}{a} \x{34}{b} ) = \psi( \x{34}{b} \x{12}{a} )$ and that
$\psi(\x{12}{a} \x{13}{b} ) = \psi( \x{13}{b} \x{12}{a} )$. The first
case is an immediate consequence of the Pure Braid Lemma \ref{pblan}
and of relation $(A1\times A1)$. Let us prove the second case, which
relies on the Pure Braid Lemma \ref{pblan} and relation $(A2)$:
$$\begin{array}{rcl}
\psi(\x{12}{a} \x{13}{b} ) &=&\uu{\y{1}{a} \ym{1}{0} }\y{2}{0} \y{1}{b}
\ym{1}{0} \ym{2}{0} \\
 &=&(\y{1}{0} )^{-2} \y{1}{a} \uu{\y{1}{0} \y{2}{0} \y{1}{b}  }\ym{1}{0}
\ym{2}{0}  \\
 &=&(\y{1}{0} )^{-2} \uu{\y{1}{a} \y{2}{b} \y{1}{0} }\, \uu{\y{2}{0}
\ym{1}{0} \ym{2}{0} } \\
 &=&(\y{1}{0} )^{-2}  \y{2}{0}  \y{1}{b} \uu{\y{2}{a}   \ym{1}{0}
\ym{2}{0} } \y{1}{0}    \\
 &=&(\y{1}{0} )^{-2}  \y{2}{0}  \y{1}{b}  \ym{1}{0}  \ym{2}{0} \uu{
\y{1}{a}   \y{1}{0} }   \\
 &=&\uu{(\y{1}{0} )^{-2} \y{2}{0}  \y{1}{b} \ym{1}{0} \ym{2}{0} (\y{1}{0}
)^2 }\uu{\y{1}{a} \ym{1}{0}}\\
 &=&\psi(\x{13}{b} )\psi(\x{12}{a} ).\\
\end{array}$$
\N $\bullet $ Relation $(St2)$: it suffices to show that $\psi(\x{12}{a}
\x{23}{b} ) =
\psi(\x{23}{b} \x{13}{ab} \x{12}{a} )$.
$$\begin{array}{rcl}
\psi(\x{12}{a} \x{23}{b} ) &=&\y{1}{a} \ym{1}{0} \uu{\y{2}{b} }\ym{2}{0} \\
&=&\y{1}{a} \uu{\ym{1}{0} \y{2}{b} \y{1}{0}}\ym{1}{0} \ym{2}{0} \\
&=&\uu{\y{1}{a} \y{2}{0} \y{1}{b}} \uu{ \ym{2}{0} \ym{1}{0} \ym{2}{0}} \\
&=&\y{2}{b} \uu{ } \y{1}{ab} \uu{ } \y{2}{a}  \ym{1}{0} \ym{2}{0}  \ym{1}{0} \\
&=&\y{2}{b} \ym{2}{0} \y{2}{0}  \y{1}{ab}\ym{1}{0} \uu{\y{1}{0}  \y{2}{a}
\ym{1}{0} \ym{2}{0} } \ym{1}{0} \\
&=&\uu{\y{2}{b} \ym{2}{0} } \uu{\y{2}{0} \y{1}{ab} \ym{1}{0}  \ym{2}{0} }
\uu{ \y{1}{a}  \ym{1}{0} } \\
&=&\psi(\x{23}{b} ) \psi(\x{13}{ab} ) \psi(\x{12}{a} ),\\
\end{array}$$
as a consequence of relation $(A2)$.

From \ref{conjan} it follows that the action of an element of $Br_n$ by
conjugation on $\Ker \pi$ depends only on its class in
${\mathcal{S}_n}$. The definition of $\psi$ on $St_n(A)$ makes clear
that it is an ${\mathcal{S}_n}$-equivariant map.

Defining $\psi$ on $Br_n$ by $\psi (y_{\alpha}) = y_{\alpha}^0 \in
Br_n(0)$ yields a group homomorphism
$$\psi : St_n(A)\rtimes Br_n \to \Ker \pi \rtimes Br_n = Br_n(A).$$

The group homomorphisms $\phi$ and $\psi$ are clearly inverse to each
other since they interchange $y_{\alpha}^a$ and $x_{\alpha}^a
y_{\alpha}$. Hence they are both isomorphisms, as asserted.  \hfill
$\square$

\begin{corollary} [Kassel-Reutenauer \cite{KR}]\label{kassel} The group
presented by generators $y_i^a$, $1\leq i\leq n-1$, $a\in A$, and
relations
\begin{align*}
 (\y{i}{0})^2&= 1 & \\
  \y{i}{a} (\y{i}{0})^{-1} \y{i}{b}&= \y{i}{a+b} & \\
  \y{i}{a} \y{j}{b} &= \y{j}{b} \y{i}{a} \qquad \text{ if } |i-j|\geq 2\\
 \y{i}{a}\y{i+1}{b}\y{i}{c}&= \y{i+1}{c}\y{i}{b+ac}\, \y{i+1}{a} & \\
\end{align*}
$a,b,c \in A$, is isomorphic to the semi-direct product $St_n(A)
\rtimes {\mathcal{S}_n}$.
\end{corollary}

Observe that when the first relation in this Corollary is deleted, the
second relation has several possible non-equivalent liftings. The one
we have chosen, $(A1)$, is what allows us to prove Theorem  \ref{homom}.

\pagebreak

\section{The parametrized braid group for $\Phi = D_n$} \label{Dn}

In this section we discuss the parametrized braid group $Br(D_n,A)$
(for a commutative ring $A$) and prove Theorem \ref{mainthm} in this
case: $Br(D_n,A)$ is isomorphic to the semi-direct product of $St( D_n,
A)$ with $Br(D_n,0)$.

\subsection{The braid group and the parametrized braid group for $D_n$}

Let $\Delta = \{\alpha_2, \alpha_{2'}, \alpha_3, \dots, \alpha_n \}$
be a fixed simple subsystem of a root system of type $D_n, n \ge 3$.
We adopt  the notation of \cite{DG} in which the simple roots on
the fork of $D_n$ are labeled $\alpha_2, \alpha_{2'}$.  The system
$D_n$ contains 2 subsystems of type $A_{n-1}$ generated by the simple
subsystems $\{\alpha_2, \alpha_3, \dots, \alpha_n \}$ and
$\{\alpha_{2'}, \alpha_3, \dots, \alpha_n \}$, and, for $n \geq 4$, a
subsystem of type $D_{n-1}$ generated by the simple subsystem
$\{\alpha_2, \alpha_{2'}, \alpha_3, \dots, \alpha_{n-1} \}$.

\setlength{\unitlength}{1cm}
\begin{picture}(10,4)
\put(2,2){\line(1,0){2}}
\put(6,2){\line(1,0){1}}
\put(1,1){\line(1,1){1}}
\put(1,3){\line(1,-1){1}}
\put(1,1.3){$2'$}
\put(1,3.3){$2$}
\put(2,2.3){$3$}
\put(3,2.3){$4$}
\put(7,2.3){$n$}
\put(0.9,0.9){$\bullet$}
\put(0.9,2.9){$\bullet$}
\put(1.9,1.9){$\bullet$}
\put(2.9,1.9){$\bullet$}
\put(3.9,1.9){$\bullet$}
\put(5.9,1.9){$\bullet$}
\put(6.9,1.9){$\bullet$}
\put(4.8,2){$\ldots$}
\end{picture}
\M

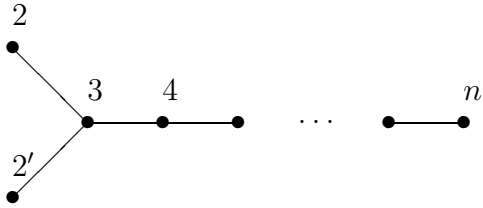
\centerline {Dynkin diagram of $D_n$}
\B

The {\it Weyl group} $W(D_n)$ is generated by the simple reflections
$\{\sigma_i = \sigma_{{\alpha}_i} \vert \alpha_i \in \Delta\}$, with
defining relations \M

$$\begin{array}{rcl}
\sigma_i^2 &=& 1 \\
(\sigma_i\sigma_j)^{m(i,j)} &=& 1\\
\end{array}$$

\M
for $i,j \in \{2,2',3, \dots, n\}$, where
$$m(i,j) = \begin{cases}  2 & \text{ if $\alpha_i,\alpha_j$ are
not connected  in the Dynkin diagram,} \\
3 & \text{ if $\alpha_i,\alpha_j$ are
 connected  in the Dynkin diagram}
\end{cases}$$

Since the only values for $m(\al,\be)$ are $1,2$ and $3$, the group
$Br(D_n, A)$ involves only relations $(A1), (A1\times A1)$ and $(A2)$.

\begin{definition}\label{brDdef} The {\it parametrized braid group of type
$D_{n}$} with parameters in the commutative ring $A$, denoted
$Br( D_n, A)$,  is generated by the elements $y_{\alpha}^a$, where
$\alpha \in \Delta$ and $a \in A$.
The relations are, for  $a, b \in A$  and $\alpha, \beta \in \Delta$
\begin{align*}
(A1)&&y_\alpha^a y_\alpha^0 y_\alpha^b& = y_\alpha^0 y_\alpha^0
y_\alpha^{a+b} & \\
(A1 \times A1)&&y_\alpha^a y_\beta^b &= y_\beta^b
y_\alpha^a &&\text{ if } m(\alpha,\beta) = 2\\
(A2)&&y_\alpha^a
y_\beta^b y_\alpha^c &= y_\beta^c y_\alpha^{b+ac} y_\beta^a &&\text{
if }m(\alpha,\beta) = 3 \text{ and } \alpha < \beta
\end{align*}
\end{definition}

Note that the simple roots in $D_n$ are ordered so that $\alpha_{2'} <
\alpha_{3}$.

\subsection{The Steinberg group of $D_n$ and the main result}

The roots of $D_n$ are $\{ \pm \epsilon_i \pm \epsilon_j \vert 1\leq
i\neq j \leq n\}$, \cite{C},\cite{H}. The Weyl group $W(D_n)\cong
(\mathbb{Z}/2)^{n-1} \rtimes {\mathcal{S}_n}$ \cite[p. 257, (X)]{Bbki}
acts on the roots by permuting the indices (action of
${\mathcal{S}_n}$) and changing the signs (action of
$(\mathbb{Z}/2)^{n-1} $). For the simple subsystem $\Delta$ we take
$\alpha_i = -\epsilon_{i-1}+\epsilon_i$ for $i= 2, \cdots , n$, and
$\alpha _{2'}= \epsilon_1 + \epsilon_2$. If $u$ and $v$ are positive
integers and $\alpha, \beta$ two roots, the linear combination
$u\alpha +v \beta$ is a root if and only if $u=1=v$, $\alpha = \pm
\epsilon_i \pm \epsilon_j , \beta \mp \epsilon_j \pm \epsilon_k $ and
$\pm \epsilon_i \pm \epsilon_k \neq 0$.  Hence in the case $\Phi =
D_n$, Definition \ref{stbg} becomes

\begin{definition}\label{stDdef}  The {\it Steinberg group of type $D_{n}$}
with parameters in the commutative ring $A$, denoted $St( D_n, A)$, is
generated by elements $x_{\alpha}^a$, where $\alpha \in D_{n}$ and $a
\in A$, subject to the relations (for $a, b \in A$ and $\alpha, \beta
\in \Phi$)
\begin{align*}
(St0)\quad &&x_\alpha^a  x_\alpha^b& = x_\alpha^{a+b} & \\
(St1)\quad &&x_\alpha^a x_\beta^b &= x_\beta^b
x_\alpha^a &\text{ if } \alpha+\beta \not\in D_n \text{ and  }
\alpha+\beta \neq 0, \\
(St2)\quad &&x_\alpha^a
x_\beta^b   &= x_\beta^b x_{\alpha+\beta}^{ab} y_\alpha^a &\text{ if }
\alpha +\beta \in  D_n.\\
\end{align*}
\end{definition}

The Weyl group $W(D_n)$ acts on $St(D_n, A)$ by $\sigma \cdot
x_\alpha^a =x_{\sigma(\alpha)}^a$ (\emph{i.e.}, in the notation of \S
1.3, one has $\eta(\delta, \gamma) = +1$) , and we can construct the
semi-direct product $St(D_n,A)\rtimes Br(D_n)$ with respect to this
action.  \M

\begin{theorem}  \label{homomD} For any commutative ring $A$ the map
$$\phi : Br(D_n,A)
\to  St(D_n,A)\rtimes Br(D_n)$$
 induced by $\phi (\y{\al}{a} ) = \x{\al}{a} \, y_{\al}{}$ is a group
isomorphism.
\end{theorem}

\begin{corollary} The group presented by generators
$y_i^a$, $i= 2', 2, 3, \cdots , n$, $a\in A$ and relations

\begin{align*}
 (\y{i}{0})^2&= 1 & \\
  \y{i}{a} (\y{i}{0})^{-1} \y{i}{b}&= \y{i}{a+b} & \\
  \y{i}{a} \y{j}{b} &= \y{j}{b} \y{i}{a}&\text{ if } |i-j|\geq 2, \text{ or
} i=2, j=2',\\
 \y{i}{a}\y{i+1}{b}\y{i}{c}&= \y{i+1}{c}\y{i}{b+ac}\, \y{i+1}{a} &
  \text{ where } i+1=3 \text{ when } i=2'\\
\end{align*}
for $a,b,c \in A$, is isomorphic to the semi-direct product $St(D_n,A)
\rtimes W(D_n)$.
\end{corollary}

\N {\it Proof of Corollary.} For each simple root $\alpha_i \in D_n$,
write $\y{i}{a}$ for $\y{\alpha_i}{a}$. \hfill \qed

\medskip
\N {\it Proof of Theorem \ref{homomD}.} The proof follows the general
pattern indicated in \S\ref{mainthmpf}.

\N {\it Step (a).} Since the relations involved in the definitions of
 $Br(D_n,A) \text{ and } St(D_n,A)$ are the same as the relations in
 the case of $A_{n-1}$, the map $\phi $ is well-defined (\emph{cf.}
 Theorem \ref{homom}).

\N {\it Step (b).} The proof of the Pure Braid Lemma in the $D_n$ case
will be given below in \S \ref{sectpbldn}.

\N {\it Step (c).} Let $\pi : Br(D_n, A) \to Br(D_n)$ be the
projection which sends each $a\in A$ to $0$ (as usual we identify $Br
(D_n,0) $ with $Br(D_n)$). We define
$$\psi : St(D_n,A) \rtimes Br(D_n) \to Br(D_n, A)\cong \Ker \pi
\rtimes Br(D_n)$$ on the first component by $\psi (\x{\alpha}{a})
=\y{\alpha}{a} \ym{\alpha}{0} \in \Ker \pi$ for $\alpha\in \Delta$.
For any $\alpha\in D_n$ there exists $\sigma \in W(D_n)$ such that
$\sigma(\alpha)\in \Delta$. Let $\tilde\sigma \in Br(D_n)$ be a
lifting of $\sigma$, and define $\psi (\x{\alpha}{a}) = \tilde
\sigma^{-1} \psi(\x{\sigma(\alpha)}{a}) \tilde\sigma \in\Ker
\pi$. This element is well-defined since it does not depend on the
lifting of $\sigma$ by the Pure Braid Lemma for $D_n$ (Lemma
\ref{PBLD}), and does not depend on the choice of $\sigma$ by Lemma
\ref{conjugation}.

In order to show that $\psi$ is a well-defined group homomorphism, it
suffices to show that the Steinberg relations are preserved. But this
is the same verification as in the $A_{n-1}$ case, (\emph{cf.}  Theorem
\ref {homom}.)  \M

The group homomorphisms $\phi$ and $\psi$ are inverse to each other
since they interchange $y_{\alpha}^a$ and $x_{\alpha}^a y_{\alpha}$. Hence they
are both isomorphisms.  \hfill $\square$ \M
\subsection{The Pure Braid Lemma for $D_n$}\label{sectpbldn}

\subsubsection{Generators for the Pure Braid Group of $D_n$} \label{gens}

In principle, the method of Reidemeister-Schreier\cite{MKS} is
available to deduce a presentation of $PBr(D_n)$ from that of
$Br(D_n)$. The details have been worked out by Digne and Gomi
\cite{DG}, although not in the specificity we need here. From their
work we can deduce that the group $PBr(D_n)$ is generated by the
elements $y_\alpha^2, \alpha \in \Delta$, together with a very small
set of their conjugates.  For example, $PBr(D_4)$ is generated by the
12 elements
$$\22^2, \22 '^2,\33^2, {}^\33\22^2, {}^\33\22 '^2, {}^{\33\22\22
'}\33^2, \44^2, {}^\44\33^2, {}^{\44\33}\22 '^2, {}^{\44\33}\22^2,
{}^{\44\33\22\22 '}\33^2, {}^{\44\33\22\22 '\33}\44^2$$ where a
prefixed exponent indicates conjugation: ${}^hg = hgh^{-1}$. Here (and
throughout) we use the simplified notations $\kk^a = y_{\alpha_k}^a$
and $\kk = y_{\alpha_k}^0$ similar to those of \S \ref{simplnot}.

\begin{proposition} \label{pbrgensdn}
For  $n \geq 4$, $PBr(D_n)$ is generated by  the elements

$\bullet \quad {\aa}_{j,i} = \jj\jmjm\dots \ipip \ii\ii \ipip
\dots\jmjm\jj, \,
n \geq j \geq i \geq 2$, and

$\bullet \quad {\bb}_{j,i} =
\jj\jmjm\dots\33\22\22'\33\dots\imim\ii\ii\imim\dots\33\22'\22\33\jmjm\jj, \,
n \geq j \geq i \geq 3$, where $i+1 = 3$ when $i = 2'$.
\end{proposition}

\N {\it Note.} Since the notation can be confusing, let us be clear
about the definition of these generators in certain special cases:

$\bullet \quad \text{When }i=j, {\aa}_{j,i} = \ii^2$.

$\bullet \quad \text{When }i=3, {\bb}_{j,3} =
\jj\jmjm\dots\33\22\22'\33\33\22'\22\dots\jmjm\jj$.

\medskip
\N {\it Proof of \eqref{pbrgensdn}.} We work in the case where $W
\text{ (in the notation of \cite{DG}) } = W(D_n)$ (in our notation;
\emph{cf.} \eqref{defnweyl}). In the proof of \cite[Corollary
2.7]{DG}, we see that $P_W = U_n \rtimes P_{W_{I_{n-1}}}$; taking $I_n
= \{\mathbf{s}_1, \mathbf{s}_2, \mathbf{s}_{2'}, \dots,
\mathbf{s}_{n-1}\}, \, n \geq 4$, as on \cite[p. 10]{DG}, we see that
their $P_W$ is equal to (our) $PBr(D_n)$ and their $P_{W_{I_{n-1}}}$ 
is equal to (our) $ PBr(D_{n-1})$.  It follows that a set of generators for
$PBr(D_n)$ can be obtained as the union of a set of generators for
$PBr(D_{n-1})$ with a set of generators for $U_n$. This sets the stage
for an inductive argument, since $D_3 = A_3$ (with $\{\alpha_2,
\alpha_3, \alpha_{2'}\} \subset D_3$ identified with $\{\alpha_1,
\alpha_2, \alpha_3\} \subset A_3$). Because $W(D_n)$ is a finite Weyl
group, it follows from \cite[Proposition 3.6]{DG}, that $U_n$ is
generated (not just normally generated) by the elements
$\mathbf{a}_{\mathbf{b},\mathbf{s}}$, and a list of these generators
in our case is given on \cite[p. 10]{DG}.

The calculations necessary to prove the Pure Braid Lemma for $D_n$ are
simpler if we replace the Digne-Gomi generators by the equivalent set
in which conjugation is replaced by {\it reflection}; that is, we
replace a generator ${}^hg = hgh^{-1}$ by $hgh'$, where if $h =
y_{i_1}\dots y_{i_k}$, $h' = y_{i_k}\dots y_{i_1}$. (We already used
this trick in the case of $A_{n-1}$.)  For $D_4$, this procedure
yields as generators of $PBr(D_4)$ the set
$$\22^2, \22'^2,\33^2, \33\22^2\33, \33\22'^2\33,
\33\22\22'\33^2\22'\22\33, \44^2, \44\33^2\44, \44\33\22^2\33\44,
\44\33\22'^2\33\44, \44\33\22\22'\33^2\22'\22\33\44,
\44\33\22\22'\33\44^2\33\22'\22\33\44$$
\N
and, more generally,  $PBr(D_n), n \geq 4$ is generated by  the elements
stated in the Proposition. \hfill \qed

\begin{lemma}[Pure
Braid Lemma for $D_n$] \label{PBLD} Let $\alpha_k \in D_n, \text{
let } y_k^a$ be a generator of $Br(D_n, A)$, and let $\omega \in
PBr(D_n)$. Then there exists $\omega '\in Br(D_n)$, independent of
$a$, such that
$$ y_k^a\  \omega = \omega '\ y_k^a.$$
Hence for any integer $k$ and any element $a \in A$,
the element $y_k^a (y_k^0)^{-1}\in Br(D_n, A)$ commutes with every
element of the pure braid group $PBr(D_n)$.
\end{lemma}

\noindent {\it Proof.} Let us show that the first assertion implies
the second one.  Let $\omega \in PBr(D_n) \subset PBr(D_n,A)$. By the
first assertion of the Lemma we have
$$y_k^a \omega = \omega ' y_k^a$$ for some $\omega ' \in Br(D_n, 0)$,
independent of $a$. Setting $a=0$ tells us that $\omega ' = y_k^0
\omega (y_k^0)^{-1}$.  Thus
$$y_k^a  (y_k^0)^{-1} \omega = \omega  y_k^a (y_k^0)^{-1} $$ as desired.
\medskip

Before beginning the proof of the first assertion, we recall some
notation introduced in \S \ref{simplnot}.  We abbreviate
$y_{\alpha_k}^a$ by $\kk^a$ and $y_{\alpha_k}^0$ by $\kk$. Whenever
there exist $\omega '$ and $\omega ''\in Br(D_n)$ such that
$\kk^a\omega = \omega ' \jj^a \omega ''$, we will write $\kk^a \omega
\sim \jj^a \omega ''$. This is {\it not} an equivalence relation, but
it is compatible with multiplication on the right by elements of
$Br(D_n)$: if $\eta \in Br(D_n)$,
$$\kk^a \omega \sim \jj^a \omega '' \Leftrightarrow \kk^a\omega = \omega '
\jj^a \omega '' \Leftrightarrow \kk^a\omega \eta = \omega ' \jj^a \omega ''
\eta
 \Leftrightarrow \kk^a \omega\eta \sim \jj^a \omega '' \eta$$

From defining relations $(A1), (A1 \times A1)$, and $(A2)$ of \S
\ref{explrelns}, we can deduce the following:

\begin{subequations}
\begin{eqnarray}
\kk^a \kk \kk & \sim & \kk^a \label{rel1} \\
 \kk^a \ \kmkm\ \kk & \sim & \kmkm^a \label{rel2} \\
  \kk^a \ \kmkm\ \kk^{-1} & \sim & \kmkm^a \label{rel3} \\
\kk^a \kpkp \kk & \sim & \kpkp^a \label{rel4} \\
 \kpkp^a \kk \kpkp  & \sim & \kk^a\label{rel5} \\
 \kpkp^a \kk \kpkp^{-1}  & \sim & \kk^a\label{rel6} \\
\kk^a \kpkp \kk^{-1} & \sim & \kpkp^a \label{rel7} \\
\kk^a \kpkp^{-1} \kk^{-1} & \sim & \kpkp^a \label{rel8} \\
\kpkp^a \kk^{-1} \kpkp^{-1}  & \sim & \kk^a \label{rel9} \\
\kk^a \kk \kmkm & \sim & \kmkm^a \kmkm \kk^{-1} \label{rel10}
\end{eqnarray}
\end{subequations}

\noindent {\it Proof, continued.} In the notation just introduced, we
must show, for every $\omega \in PBr(D_n)$, that $\kk^a \omega \sim
\kk^a$, and it suffices to show this when $\omega$ is one of the
generators ${\aa}_{j,i}$ or ${\bb}_{j,i}$ of \S \ref{gens}. That is,
we must show

\begin{subequations}
\begin{eqnarray}
\quad \kk^a {\aa}_{j,i} &\sim& \kk^a  \quad n \geq j \geq i \geq 2,\quad 1
\leq k \leq n \label{reln1} \\
\quad \kk^a {\bb}_{j,i} &\sim& \kk^a  \quad  n \geq j \geq i \geq
3,\quad 1 \leq k \leq n \label{reln2}
\end{eqnarray}
\end{subequations}

The proofs of \eqref{reln1} for $i \geq 3$ are exactly the same as the
corresponding proofs for $A_{n-1}$ (see \S 2); the additional case
$i=2'$ presents no new issues.  Thus we shall concentrate on proving
\eqref{reln2}; the proof proceeds by induction on $n$.

The case $n=3$ is the case of the root system $D_3 = A_3$, which is
part of the Pure Braid Lemma \ref{pblan} for $A_{n-1}$.  Hence we may
assume $n \geq 4$, and that \eqref{reln2} holds whenever $j, k \leq
n-1$.  That is, we must prove \eqref{reln2} in these cases:
\begin{equation*}
k=n, \, j \leq n-1; \quad
k \leq n-1, \, j=n; \quad
k = n, \, j=n
\end{equation*}
which further subdivide into the cases
\begin{eqnarray}
&k=n ,& j \leq n-2 \label{eqn1} \\
&k=n ,&  j =  n-1 \label{eqn2}  \\
&k \leq n-2 ,&  j=n \label{eqn3}  \\
&k =  n-1 ,&  j=n \label{eqn4}  \\
&k = n ,&  j=n \label{eqn5}
\end{eqnarray}

\noindent $\bullet$ {\it Case \eqref{eqn1} $k=n$ and $ j \leq n-2 $.}
Since $i \leq j  \leq n-2$, it follows from relation
$(A1 \times A1)$ that $\nn^a$ commutes with every generator which
occurs in the expression for ${\bb}_{j,i}$; hence \begin{eqnarray*}
\nn^a{\bb}_{j,i} &=&
\nn^a\jj\jmjm\dots\33\22\22'\33\dots\imim\ii\ii\imim\dots\33\22'\22\33\jmjm\jj
\\ &=& {\bb}_{j,i} \nn^a \\ &\sim & \nn^a
\end{eqnarray*} as desired.

\medskip \noindent $\bullet$ {\it Case \eqref{eqn2}  $k=n$ and $ j =  n-1$.} \\

\noindent
If $i \leq n-2$, then
\begin{eqnarray*}
\nn^a{\bb}_{n-1,i} &=&
\nn^a\nmnm\nmmnmm\dots\33\22\22'\33\dots\imim\ii\ii\imim\dots\33\22'\22\33\dots\nmmnmm\nmnm \\ &=& \underbrace{\nn^a
\nmnm\nn^{-1}}\nn\nmmnmm\dots\33\22\22'\33\dots\imim\ii\ii\imim\dots\33\22'\22\33\dots\nmmnmm\nmnm \\ &\sim & \nmnm^a
\nn\nmmnmm\dots\33\22\22'\33\dots\imim\ii\ii\imim\dots\33\22'\22\33\dots\nmmnmm\nmnm \\ &=& \underbrace{ \nmnm^a
\nmmnmm\dots\33\22\22'\33\dots\imim\ii\ii\imim\dots\33\22'\22\33\dots\nmmnmm}\nn
\nmnm \\ &\sim & \nmnm^a \nn\nmnm \\ &\sim & \nn^a
\end{eqnarray*}
 as desired.

The case $k=n, i=j=n-1$, is considerably more complicated.  We first
prove some preliminary lemmas.

\begin{lemma}\label{twin}
$$\nn^a\nmnm\nmmnmm\dots\33\22\22' \sim
 \22^a\33\22'\44\dots\nmmnmm\nmnm\nn $$
\end{lemma}
\noindent {\it Proof.}
\begin{eqnarray*}
\nn^a\nmnm\nmmnmm\dots\33\22\22' & = &
\nn^a\nmnm\nn^{-1}\nn\nmmnmm\dots\33\22\22' \\
& \sim & \nmnm^a\nn\nmmnmm\dots\33\22\22' \\
& = & \nmnm^a\nmmnmm\dots\33\22\22'\nn \\
& & \vdots\\
& \sim & \33^a\22\22'\44\dots\nmmnmm\nmnm\nn \\
& = &  \33^a\22\33^{-1}\33\22'\44\dots\nmmnmm\nmnm\nn \\
& \sim & \22^a\33\22'\44\dots\nmmnmm\nmnm\nn \\
\end{eqnarray*}
\hfill $\square$

\begin{lemma} \label{twine}
$$\33\22'\44\33\55\44\dots\nmmmnmmm\nmnm\nmmnmm\nn\nn =
\33\44\55\dots\nmnm\nn\nn\22'\33\44\dots\nmmnmm $$
\end{lemma}
\noindent {\it Proof.}
\begin{eqnarray*}
\33\22'\44\33\55\44\dots\nmmmnmmm\nmnm\nmmnmm\nn\nn
& = &\33\22'\44\33\55\44\dots\nmmmnmmm\nmnm\nn\nn\nmmnmm \\
& = &\33\22'\44\33\55\44\dots\nmnm\nn\nn\nmmmnmmm\nmmnmm \\
& & \vdots\\
& = &\33\44\55\dots\nmnm\nn\nn\22'\33\44\dots\nmmnmm
\end{eqnarray*}
\hfill $\square$

We now complete the case $k=n, i=j=n-1$.
\begin{eqnarray*}
\nn^a{\bb}_{n-1,n-1}
&=&\nn^a\nmnm\nmmnmm\dots\33\22\22'\33\dots\nmmnmm\nmnm\nmnm\nmmnmm\dots\\
&&\dots\33\22'\22\33\dots\nmmnmm\nmnm \\
&\sim&\22^a\33\22'\44\dots\nmmnmm\nmnm\nn
\33\dots\nmmnmm\nmnm\nmnm\nmmnmm\dots \\
&&\dots\33\22'\22\33\dots\nmmnmm\nmnm \quad \text{(by Lemma
\ref{twin})} \\
&=&\22^a\33\22'\44\dots\nmmnmm\nmnm\33\dots\nmmnmm\underbrace{\nn
\nmnm\nn^{-1}} \nn \nmnm\nmmnmm\dots \\
&&\dots\33\22'\22\33\dots\nmmnmm\nmnm \\
&=&\22^a\33\22'\44\dots\nmmnmm\underbrace{\nmnm\33\dots}\nmmnmm\nmnm^{-1}\nn\nmnm
\nn \nmnm\nmmnmm\dots \\ & &\dots\33\22'\22\33\dots\nmmnmm\nmnm \\
&=&\22^a\33\22'\44\dots\nmmnmm\33\dots\underbrace{\nmnm\nmmnmm\nmnm^{-1}}\nn\nmnm
\nn \nmnm\nmmnmm\dots \\ & &\dots\33\22'\22\33\dots\nmmnmm\nmnm \\
&=&\22^a\33\22'\44\dots\nmmnmm\33\dots\underbrace{\nmmnmm^{-1}\nmnm\nmmnmm}\nn\nmnm \nn \nmnm\nmmnmm\dots \\ & &\dots\33\22'\22\33\dots\nmmnmm\nmnm \\
&&\vdots\\
&=&\22^a\33\22'\underbrace{\44\33\44^{-1}}\55\44\dots\nmnm\nmmnmm\nn\nmnm
\nn \nmnm\nmmnmm\dots \\ & &\dots\33\22'\22\33\dots\nmmnmm\nmnm \\
&=&\22^a\underbrace{\33\22'\33^{-1}}\44\33\55\44\dots\nmnm\nmmnmm\nn\nmnm
\nn \nmnm\nmmnmm\dots \\ & &\dots\33\22'\22\33\dots\nmmnmm\nmnm \\
&=&\22^a\underbrace{\22'^{-1}\33\22'}\44\33\55\44\dots\nmnm\nmmnmm\nn\nmnm
\nn \nmnm\nmmnmm\dots \\ & &\dots\33\22'\22\33\dots\nmmnmm\nmnm \\
&\sim&\22^a\33\22'\44\33\55\44\dots\nmnm\nmmnmm\nn\underbrace{\nmnm
\nn \nmnm}\nmmnmm\dots \\ & &\dots\33\22'\22\33\dots\nmmnmm\nmnm \\
&=&\22^a\33\22'\44\33\55\44\dots\nmnm\nmmnmm\nn\underbrace{\nn
\nmnm\nn}\nmmnmm\dots \\ & &\dots\33\22'\22\33\dots\nmmnmm\nmnm \\
&=&\22^a\33\44\55\dots\nmnm\nn\nn\22'\33\44\dots\nmmnmm
\nmnm\nn\nmmnmm\dots \\ & &\dots\33\22'\22\33\dots\nmmnmm\nmnm \quad
\text{(by Lemma \ref{twine})}
\end{eqnarray*}

We now manipulate part of this expression so that we can apply induction.
\begin{eqnarray*}
\22^a\33\44\55\dots\nmnm\nn\nn\22'
& =
&\22^a\33\44\55\dots\nmnm\nn\nn\underbrace{\nmnm\dots\55\44\33\33^{-1}\44^{-1}\55^{-1}\dots\nmnm^{-1}}\22 ' \\
& \sim &\22^a\33^{-1}\44^{-1}\55^{-1}\dots\nmnm^{-1}\22 '
\end{eqnarray*}
(by the case of $A_{n-1} = \{\alpha_2, \alpha_3, \dots, \alpha_n\}$).
Hence
\begin{eqnarray*}
\nn^a{\bb}_{n-1,n-1}
&=&\22^a\33^{-1}\44^{-1}\55^{-1}\dots\nmnm^{-1}\22'\33\44\dots\nmmnmm
\nmnm\nn\nmmnmm\nmmmnmmm\dots \\
& &\dots\33\22'\22\33\dots\nmmnmm\nmnm \\
&=&\22^a\33^{-1}\44^{-1}\55^{-1}\dots\nmnm^{-1}\22'\33\44\dots\nmmnmm
\nmnm\nmmnmm\nn\nmmmnmmm\dots \\
& &\dots\33\22'\22\33\dots\nmmnmm\nmnm \\
&=&\22^a\33^{-1}\44^{-1}\55^{-1}\dots\nmnm^{-1}\22'\33\44\dots
\nmnm\nmmnmm\nmnm\nn\nmmmnmmm\dots \\
& &\dots\33\22'\22\33\dots\nmmnmm\nmnm \\
&=&\22^a\33^{-1}\44^{-1}\55^{-1}\dots\nmmnmm^{-1}\22'\33\44\dots
\nmmnmm\nmnm\nn\nmmmnmmm\dots \\
& &\dots\33\22'\22\33\dots\nmmnmm\nmnm \\
& &\vdots\\
&=&\22^a\33^{-1}\22'\33\44\dots\nmmnmm\nmnm\nn\22'\22\33\dots\nmmnmm\nmnm \\
&=&\22^a\33^{-1}\22'\33\22'\44\dots\nmmnmm\nmnm\nn\22\33\dots\nmmnmm\nmnm \\
&=&\22^a\33^{-1}\33\22'\33\44\dots\nmmnmm\nmnm\nn\22\33\dots\nmmnmm\nmnm \\
&=&\22^a\22'\33\44\dots\nmmnmm\nmnm\nn\22\33\dots\nmmnmm\nmnm \\
&\sim&\22^a\33\44\dots\nmmnmm\nmnm\nn\22\33\dots\nmmnmm\nmnm \\
&=&\22^a\33\22\44\dots\nmmnmm\nmnm\nn\33\dots\nmmnmm\nmnm \\
&\sim&\33^a\44\dots\nmmnmm\nmnm\nn\33\dots\nmmnmm\nmnm \\
& &\vdots\\
&\sim&\nmnm^a\nn\nmnm \\
&\sim&\nn^a\\
\end{eqnarray*} as desired.\hfill $\square$

\medskip
\noindent
$\bullet$ {\it Case \eqref{eqn3} $k \leq n-2$ and $j=n$.}
\begin{eqnarray*}\
\kk^a{\bb}_{n,i} &=&
\kk^a\nn\nmnm\dots\33\22\22'\33\dots\imim\ii\ii\imim\dots\33\22'\22\33
\dots\nmnm
\nn \\
&\sim&
\underbrace{\kk^a\kpkp\dots\33\22\22'\33\dots\imim\ii\ii\imim\dots\33\22'\22\33\dots\kpkp}\kppkpp\dots\nmnm\nn \\
 \\
&\sim& \kk^a\kppkpp\dots\nn \quad
\text{(by induction)} \\
&\sim& \kk^a
\end{eqnarray*}as desired.
\smallskip

\noindent
$\bullet$ {\it Case  \eqref{eqn4} $k=n-1$ and $j=n$.}
\begin{align*}
\nmnm^a{\bb}_{n,i} &=
\nmnm^a\nn\nmnm\dots\33\22\22'\33\dots\imim\ii\ii\imim\dots\33\22'\22\33\dots\nmnm\nn
\\ &=
\underbrace{\nmnm^a\nn\nmnm}\dots\33\22\22'\33\dots\imim\ii\ii\imim\dots\33\22'\
22\33\dots\nmnm\nn \\ &\sim
\nn^a\nmmnmm\dots\33\22\22'\33\dots\imim\ii\ii\imim\dots\33\22'\22\33\dots\nmnm\nn
\tag{$\dag$}
\end{align*}

Suppose first that $i \leq n-2$.  Then
\begin{eqnarray*}
(\dag) &=&\nn^a\nmmnmm\dots\33\22\22
'\33\dots\imim\ii\ii\imim\dots\33\22 '\22\33\dots\nmnm\nn
\\ &=&
\nmmnmm\dots\33\22\22 '\33\dots\imim\ii\ii\imim\dots\33\22 '\22\33\dots\nmmnmm\nn^a\nmnm\nn
\\ &\sim&\nn^a\nmnm\nn
\\ &\sim&\nmnm^a
\end{eqnarray*} as desired.
\smallskip

\noindent
If $i = n-1$
\begin{eqnarray*}
(\dag)
&=& \nn^a\nmmnmm\dots\33\22\22'\33\dots\nmmnmm\nmnm\nmnm\nmmnmm\dots \\
& &\dots\33\22'\22\33\dots\nmnm\nn \\
&\sim& \underbrace{\nn^a\nmnm\nmnm}\nmmnmm\dots \33\22'\22\33\dots\nmnm\nn \\
&\sim& \nn^a\nmmnmm\dots \33\22'\22\33\dots\nmnm\nn  \quad
\text{(by the case $A_2 = \{\alpha_{n-1},\alpha_n\}$)}\\
&\sim& \nn^a\nmnm\nn \quad
\text{(by relation $(A1 \times A1)$}\\
&\sim&\nmnm^a
\end{eqnarray*} as desired.
\smallskip

If $i = n$
\begin{eqnarray*}
(\dag)
&=&\nn^a\nmmnmm\dots\33\22\22'\33\dots\nmnm\nn\nn\nmnm\dots
\33\22'\22\33\dots\nmnm\nn \\
&\sim&\nn^a\nmnm\nn\nn\nmnm\nmmnmm\dots \33\22'\22\33\dots\nmnm\nn \\
&\sim&\nmnm^a\nn\nmnm\nmmnmm\dots \33\22'\22\33\dots\nmnm\nn \\
&\sim&\nn^a\nmmnmm\dots \33\22'\22\33\dots\nmnm\nn \\
&\sim&\nn^a\nmnm\nn \\
&\sim&\nmnm^a
\end{eqnarray*}as desired.

\medskip
\noindent
$\bullet$ {\it Case \eqref{eqn5} $k = n$ and $j=n$.}
\begin{align*}
\nn^a{\bb}_{n,i}
&=
\nn^a\nn\nmnm\nmmnmm\dots\33\22\22'\33\dots\imim\ii\ii\imim\dots\33\22'\22\33\dots\nmnm\nn
\\
&\sim
\nmnm^a\nmnm\nn^{-1}\nmmnmm\dots\33\22\22'\33\dots\imim\ii\ii\imim\dots\33\22'\22\33\dots\nmnm\nn \tag{$\dag$}
\\
\end{align*}

If $i \leq n-2$, then
\begin{align*}
(\dag)
&=\nmnm^a\nmnm\nmmnmm\dots\33\22\22'\33\dots\imim\ii\ii\imim\dots\33\22'\22\33\dots\underbrace{\nn^{-1}\nmnm\nn}
\\
&=\underbrace{\nmnm^a\nmnm\nmmnmm\dots\33\22\22'\33\dots\imim\ii\ii\imim\dots\33\22'\22\33\dots\nmnm}\nn\nmnm^{-1}
\\ &\sim \nmnm^a\nn\nmnm^{-1} \quad \text{(by induction)} \\ &\sim
\nn^a \\
\end{align*} as desired.
\medskip

If $i = n-1$, then
\begin{eqnarray*}
(\dag)
&=&\nmnm^a\nmnm\nmmnmm\dots\33\22\22'\33\dots\nmmnmm\nn^{-1}\nmnm\nmnm\nmmnmm\dots\33\22'\22\33\dots\nmnm\nn
\\ & & \vdots \\
&\sim&\33^a\33\22\22'\44^{-1}\33\55^{-1}\44\dots\nmnm^{-1}\nmmnmm\nn^{-1}\nmnm\nmnm\nmmnmm\dots\33\22'\22\33\dots\nmnm\nn
\\
&=&\33^a\33\22\22'\44^{-1}\33\55^{-1}\44\dots\nmnm^{-1}\nmmnmm\underbrace{\nn^{-1}\nmnm\nn}\nn^{-1}\nmnm\nmmnmm\dots\\
& &\dots\33\22'\22\33\dots\nmnm\nn \\ &=&
\33^a\33\22\22'\44^{-1}\33\55^{-1}\44\dots\underbrace{\nmnm^{-1}\nmmnmm\nmnm}\nn
\underbrace{\nmnm^{-1}\nn^{-1}\nmnm}\nmmnmm\dots\\ &
&\dots\33\22'\22\33\dots\nmnm\nn \\ &=&
\33^a\33\22\22'\44^{-1}\33\55^{-1}\44\dots\nmmmnmmm\underbrace{\nmmnmm\nmnm\nmmnmm^{-1}}\nn\underbrace{\nn\nmnm^{-1}}\nn^{-1}\nmmnmm\dots\\
& &\dots\33\22'\22\33\dots\nmnm\nn \\ &=&
\33^a\33\22\22'\44^{-1}\33\55^{-1}\44\dots\nmmmnmmm\nmmnmm\nmnm\nn\nn\nmmnmm^{-1
}\nmnm^{-1}\nn^{-1}\nmmnmm\dots\\ & &\dots\33\22'\22\33\dots\nmnm\nn
\\ & &\vdots \\ &=&
\33^a\33\22\22'\underbrace{\44^{-1}\33\44}\dots\nmmnmm\nmnm\nn\nn\44^{-1}\55^{-1
}\dots \nmnm^{-1}\nn^{-1}\nmmnmm\dots\33\22'\22\33\dots\nmnm\nn \\
&\sim&
\22^a\22\underbrace{\33^{-1}\22'\33}\44\33^{-1}\55\dots\nmnm\nn\nn\44^{-1}\55^{-
1}\dots \nmnm^{-1}\nn^{-1}\nmmnmm\dots\33\22'\22\33\dots\nmnm\nn \\
&=&
\22^a\22\22'\33\22'^{-1}\44\dots\nmnm\nn\nn\33^{-1}\44^{-1}\55^{-1}\dots
\nmnm^{-1}\nn^{-1}\nmmnmm\dots\33\22'\22\33\dots\nmnm\nn \\ &=&
\22^a\22\22'\33\44\dots\nmnm\nn\nn\22'^{-1}\33^{-1}\44^{-1}\55^{-1}\dots
\nmnm^{-1}\nn^{-1}\nmmnmm\dots\33\22'\22\33\dots\nmnm\nn \\ &\sim&
\22^a\22\33\44\dots\nmnm\nn\nn\22'^{-1}\33^{-1}\44^{-1}\55^{-1}\dots
\nmnm^{-1}\nn^{-1}\nmmnmm\dots\33\22'\22\33\dots\nmnm\nn \\ &=&
\22^a\22\33\44\dots\nmnm\nn\nn\underbrace{\nmnm\dots\44\33\22\22^{-1}\33^{-1}\44
^{-1}\dots\nmnm^{-1}}\22'^{-1}\33^{-1}\44^{-1}\55^{-1}\dots \\ &
&\dots\nmnm^{-1}\nn^{-1}\nmmnmm\dots\33\22'\22\33\dots\nmnm\nn \\
&\sim&
\22^a\22^{-1}\33^{-1}\44^{-1}\dots\nmnm^{-1}\22'^{-1}\33^{-1}\44^{-1}\55^{-1}\do
ts \nmnm^{-1}\nn^{-1}\nmmnmm\dots\33\22'\22\33\dots\nmnm\nn \\ &
&\text{(by the case of $A_{n-1} = \{\alpha_2, \alpha_3, \dots,
\alpha_n\}$)} \\ &\sim&
\22^a\22^{-1}\33^{-1}\44^{-1}\dots\nmnm^{-1}\22'^{-1}\33^{-1}\44^{-1}\55^{-1}\dots
\nmmmnmmm^{-1}\underbrace{\nmmnmm^{-1}\nmnm^{-1}\nmmnmm}\nn^{-1}\nmmmnmmm\dots
\\ & &\dots\33\22 '\22\33\dots\nmnm\nn \\ &=&
\22^a\22^{-1}\33^{-1}\44^{-1}\dots\nmnm^{-1}\22
'^{-1}\33^{-1}\44^{-1}\55^{-1}\do ts
\nmmmnmmm^{-1}\underbrace{\nmnm\nmmnmm^{-1}\nmnm^{-1}}\nn^{-1}\nmmmnmmm\dots
\\ & &\dots\33\22'\22\33\dots\nmnm\nn \\ &=&
\22^a\22^{-1}\33^{-1}\44^{-1}\dots\nmnm^{-1}\nmnm\22'^{-1}\33^{-1}\44^{-1}\55^{-
1}\dots \nmmmnmmm^{-1}\nmmnmm^{-1}\nmnm^{-1}\nn^{-1}\nmmmnmmm\dots \\
& &\dots\33\22'\22\33\dots\nmnm\nn \\ &=&
\22^a\22^{-1}\33^{-1}\44^{-1}\dots\nmmnmm\22'^{-1}\33^{-1}\44^{-1}\55^{-1}\dots
\nmmmnmmm^{-1}\nmmnmm^{-1}\nmnm^{-1}\nn^{-1}\nmmmnmmm\dots \\ &
&\dots\33\22'\22\33\dots\nmnm\nn \\ & &\vdots \\ & =
&\22^a\22^{-1}\33^{-1}\22'^{-1}\33^{-1}\44^{-1}\dots\nn^{-1}
\22'\22\33\dots\nmnm\nn \\ & =
&\22^a\22^{-1}\33^{-1}\underbrace{\22'^{-1}\33^{-1}\22'}\44^{-1}\dots\nn^{-1}
\22\33\dots\nmnm\nn
\\ & = &\22^a\22^{-1}\33^{-1}\33\22
'^{-1}\33^{-1}\44^{-1}\dots\nn^{-1} \22\33\dots\nmnm\nn \\ &\sim
&\22^a\22^{-1}\33^{-1}\44^{-1}\dots\nn^{-1} \22\33\dots\nmnm\nn \\ &=
&\22^a\underbrace{\22^{-1}\33^{-1}\22}\44^{-1}\dots\nn^{-1}
\33\dots\nmnm\nn \\ &=
&\underbrace{\22^a\33\22^{-1}}\33^{-1}\44^{-1}\dots\nn^{-1}
\33\dots\nmnm\nn \\ &\sim&\33^a\33^{-1}\44^{-1}\dots\nn^{-1}
\33\dots\nmnm\nn \quad \text{ (by \eqref{rel7})}\\ & & \vdots \\
&\sim&\nmnm^a\nmnm^{-1}\underbrace{\nn^{-1}\nmnm\nn}\quad \text{(by
\eqref{rel7})}\\ &=&\nmnm^a\nmnm^{-1}\nmnm\nn\nmnm^{-1} \\
&=&\nmnm^a\nn\nmnm^{-1} \\ &\sim&\nn^a \\
\end{eqnarray*} as desired.

\medskip

If $i = n  $, then
\begin{eqnarray*}
(\dag)
&=&\nmnm^a\nmnm\nn^{-1}\nmmnmm\dots\33\22\22'\33\dots\nmnm\nn\nn\nmnm\dots\33\22
'\22\33\dots\nmnm\nn
\\
&=&\nmnm^a\nmnm\nmmnmm\dots\33\22\22'\33\dots\underbrace{\nn^{-1}\nmnm\nn}\nn\nmnm\dots\33\22'\22\33\dots\nmnm\nn
\\
&=&\nmnm^a\nmnm\nmmnmm\dots\33\22\22 '\33\dots\nmmnmm\nmnm\nn\underbrace{\nmnm^{-
1}\nn\nmnm}\dots\\
& &\dots\33\22'\22\33\dots\nmnm\nn \\
&=&\nmnm^a\nmnm\nmmnmm\dots\33\22\22'\33\dots\nmmnmm\nmnm\nn\nn\nmnm\nn^{-1}\nmmnmm\dots
\\ & &\dots\33\22'\22\33\dots\nmnm\nn \\
&=&\nmnm^a\nmnm\nmmnmm\dots\33\22\22'\33\dots\nmmnmm\nmnm\nn\nn\nmnm\nmmnmm\dots
\\ & &\dots\33\22'\22\33\dots\nn^{-1}\nmnm\nn \\
&=&\nmnm^a\nmnm\nmmnmm\dots\33\22\22'\33\dots\nmmnmm\nmnm\nn\nn\nmnm\nmmnmm\dots
\\ & &\dots\33\22'\22\33\dots\nmnm\nn\nmnm^{-1} \\
&\sim&\nmmnmm^a\nmmnmm\nmnm^{-1}\nmmmnmmm\dots\33\22\22'\33\dots\nmmnmm\nmnm\nn\nn\nmnm\nmmnmm\dots
\\ & &\dots\33\22'\22\33\dots\nmnm\nn\nmnm^{-1}  \quad \text{ (by
\eqref{rel7})}\\
&=&\nmmnmm^a\nmmnmm\nmmmnmmm\dots\33\22\22'\33\dots\nmnm^{-1}\nmmnmm\nmnm\nn\nn\nmnm\nmmnmm\dots
\\ & &\dots\33\22'\22\33\dots\nmnm\nn\nmnm^{-1} \\
&=&\nmmnmm^a\nmmnmm\nmmmnmmm\dots\33\22\22'\33\dots\nmmnmm\nmnm\nmmnmm^{-1}\nn\nn\nmnm\nmmnmm\dots
\\ & &\dots\33\22'\22\33\dots\nmnm\nn\nmnm^{-1} \\
&=&\nmmnmm^a\nmmnmm\nmmmnmmm\dots\33\22\22'\33\dots\nmmnmm\nmnm\nn\nn\nmmnmm^{-1}\nmnm\nmmnmm\dots
\\ & &\dots\33\22'\22\33\dots\nmnm\nn\nmnm^{-1} \\
&=&\nmmnmm^a\nmmnmm\nmmmnmmm\dots\33\22\22'\33\dots\nmmnmm\nmnm\nn\nn
\nmnm\nmmnmm\nmnm^{-1}\dots \\ &
&\dots\33\22'\22\33\dots\nmnm\nn\nmnm^{-1} \\
&=&\nmmnmm^a\nmmnmm\nmmmnmmm\dots\33\22\22'\33\dots\nmmnmm\nmnm\nn\nn
\nmnm\nmmnmm\dots \\ &
&\dots\33\22'\22\33\dots\nmmmnmmm\nmnm^{-1}\nmmnmm\nmnm\nn\nmnm^{-1}
\\
&=&\nmmnmm^a\nmmnmm\nmmmnmmm\dots\33\22\22'\33\dots\nmmnmm\nmnm\nn\nn
\nmnm\nmmnmm\dots \\ &
&\dots\33\22'\22\33\dots\nmmmnmmm\nmmnmm\nmnm\nmmnmm^{-1}\nn\nmnm^{-1}
\\
&=&\nmmnmm^a\nmmnmm\nmmmnmmm\dots\33\22\22'\33\dots\nmmnmm\nmnm\nn\nn
\nmnm\nmmnmm\dots \\ &
&\dots\33\22'\22\33\dots\nmmmnmmm\nmmnmm\nmnm\nn\nmmnmm^{-1}\nmnm^{-1}
\\ & & \vdots \\
&\sim&\22^a\22\22'\33\dots\nmnm\nn\nn\nmnm\dots\33\22'\22\33\dots\nmnm\nn\22^{-1
}\33^{-1}\dots\nmnm^{-1} \quad \text{ (by \eqref{rel7})}
\\
&\sim&\underbrace{\22^a\22\33\dots\nmnm\nn\nn\nmnm\dots\33\22}\22'\33\dots\nmnm\
nn\22^{-1}\33^{-1}\dots\nmnm^{-1}
\\
&\sim&\22^a\22'\33\dots\nmnm\nn\22^{-1}\33^{-1}\dots\nmnm^{-1}\quad \text{
(by the case of $A_{n-1} = \{\alpha_2, \alpha_3, \dots,
\alpha_n\}$)} \\
&\sim&\22^a\33\dots\nmnm\nn\22^{-1}\33^{-1}\dots\nmnm^{-1} \\
&=&\22^a\33\22^{-1}\dots\nmnm\nn\33^{-1}\dots\nmnm^{-1} \\
&=&\33^a\44\dots\nmnm\nn\33^{-1}\dots\nmnm^{-1} \\
\\ & & \vdots \\
&=&\nmnm^a\nn\nmnm^{-1} \\
&=&\nn^a
\end{eqnarray*} as desired.

\pagebreak

\B
 Institut de Recherche Math\'ematique Avanc\'ee,

    CNRS et Universit\'e Louis Pasteur

    7 rue R. Descartes,

    67084 Strasbourg Cedex, France

    Courriel: loday@math.u-strasbg.fr

\B
\N

and

\B
Department of Mathematics

Northwestern University

2033 Sheridan Road

Evanston IL 60208--2730  USA

e-mail: mike@math.northwestern.edu
\end{document}

%% file: LodayStein0.tex
\usepackage[fleqn]{amsmath}
\usepackage{amsfonts}
\usepackage{amssymb}
\usepackage{amsthm}
\def\11{\mathbf{1}}
\def\22{\mathbf{2}}
\def\33{\mathbf{3}}
\def\44{\mathbf{4}}
\def\55{\mathbf{5}}
\def\kk{\mathbf{k}}
\def\nn{\mathbf{n}}
\def\aa{\mathbf{a}}
\def\bb{\mathbf{b}}
\def\jj{\mathbf{j}}
\def\ii{\mathbf{i}}
\def\ipip{\mathbf{(i+1)}}
\def\jpjp{\mathbf{(j+1)}}
\def\jmjm{\mathbf{(j-1)}}
\def\imim{\mathbf{(i-1)}}
\def\kmkm{\mathbf{(k-1)}}
\def\kmmkmm{\mathbf{(k-2)}}
\def\kpkp{\mathbf{(k+1)}}
\def\kppkpp{\mathbf{(k+2)}}
\def\nmnm{\mathbf{(n-1)}}
\def\nmmnmm{\mathbf{(n-2)}}
\def\nmmmnmmm{\mathbf{(n-3)}}
\allowdisplaybreaks[3]

\def\al{\alpha}
\def\be{\beta}

\def\DD{\Delta}

\def\ss{\sigma}

\def\oo{\omega}

\def\M{\medskip}
\def\B{\bigskip}
\def\BB{\bigskip\bigskip}
\def\N{\noindent}

\def\Ker{\mathop{\rm Ker\,}\nolimits}

\newcommand{\y}[2]{y_{#1}^{#2}}
\newcommand{\ym}[2]{(y_{#1}^{#2})^{-1}}

\newcommand{\uu}{\underbrace }
\newcommand{\x}[2]{x_{#1}^{#2}}

\def\mm{\mathop{^{-1}}\nolimits}

\newtheorem{theorem}{Theorem}[subsection]
\newtheorem{proposition}[theorem]{Proposition}
\newtheorem{lemma}[theorem]{Lemma}
\newtheorem{corollary}[theorem]{Corollary}
\newtheorem{definition}[theorem]{Definition}